\documentclass[oneside]{amsart}

\usepackage[letterpaper, total={6.5in, 8.5in}]{geometry}
\usepackage{amsmath}
\usepackage{amsfonts}
\usepackage{amssymb}
\usepackage[dvipsnames]{xcolor}
\usepackage{mathrsfs}
\usepackage{enumerate}
\usepackage{amsthm}
\usepackage{xcolor}
\usepackage{color}
\usepackage[pdfdisplaydoctitle=true,
            colorlinks=true,
            urlcolor=blue,
            citecolor=blue,
            linkcolor=blue,
            pdfstartview=FitH,
            pdfpagemode= UseNone,
            bookmarksnumbered=true]{hyperref}
\usepackage{tikz}
\usetikzlibrary{decorations.pathreplacing}
\usepackage{mathrsfs}
\usepackage{pgfplots}
\usepgfplotslibrary{fillbetween}
\pgfplotsset{compat=1.18}
\usepackage[normalem]{ulem}
\usepackage{thmtools}

\definecolor{darkgreen}{rgb}{0,0.5,0}
\definecolor{darkred}{rgb}{0.7, 0, 0}

\makeatletter
\def\@tocline#1#2#3#4#5#6#7{\relax
  \ifnum #1>\c@tocdepth 
  \else
    \par \addpenalty\@secpenalty\addvspace{#2}%
    \begingroup \hyphenpenalty\@M
    \@ifempty{#4}{%
      \@tempdima\csname r@tocindent\number#1\endcsname\relax
    }{%
      \@tempdima#4\relax
    }%
    \parindent\z@ \leftskip#3\relax \advance\leftskip\@tempdima\relax
    \rightskip\@pnumwidth plus4em \parfillskip-\@pnumwidth
    #5\leavevmode\hskip-\@tempdima
      \ifcase #1
       \or\or \hskip 1em \or \hskip 2em \else \hskip 3em \fi%
      #6\nobreak\relax
    \hfill\hbox to\@pnumwidth{\@tocpagenum{#7}}\par
    \nobreak
    \endgroup
  \fi}
\makeatother

\newtheorem{theorem}{Theorem}[section]
\newtheorem{lemma}[theorem]{Lemma}
\newtheorem{proposition}[theorem]{Proposition}

\theoremstyle{definition}

\newtheorem*{theorem*}{Theorem}
\theoremstyle{remark}
\newtheorem{remark}[theorem]{Remark}

\numberwithin{equation}{section}
\newcommand{\R}{\mathbb{R}}
\newcommand{\T}{\mathbb{T}}

\newcommand{\Z}{\mathbb{Z}}
\newcommand{\C}{\mathbb{C}}

\newcommand{\D}{\mathscr{D}}
\newcommand{\Dmu}{\mathscr{D}_{\mu}}
\newcommand{\Dex}{\mathscr{D}_{\mu,\text{ex}}}

\newcommand{\Amu}{\mathscr{A}_{\mu}}
\newcommand{\Swi}{\mathscr{A}_{\mu, 0}}

\newcommand{\symp}{\nabla^\perp}
\DeclareMathOperator{\curl}{curl}
\renewcommand{\div}{\mathrm{div}}
\DeclareMathOperator{\Ad}{Ad}
\DeclareMathOperator{\ad}{ad}

\newcommand{\abs}[1]{\lvert#1\rvert}
\newcommand{\norm}[1]{\left\lVert#1\right\rVert}
\newcommand{\ip}[1]{\langle #1 \rangle}

\newcommand{\bigip}[1]{\big\langle #1 \big\rangle}

\newcommand{\Adstar}[1]{\Ad_{#1}^\star}
\newcommand{\adstar}[1]{\ad_{#1}^\star}

\newcommand{\radstar}[1]{\ad_{#1}^{\star_r}}

\author[1]{James Benn}
\author[2]{Patrick Heslin}
\author[3]{Leandro Lichtenfelz}
\author[4]{Gerard Misio{\l}ek}

\address{J. Benn: Melbourne, Australia}
\email{jbenn2@alumni.nd.edu}
\address{P. Heslin: Max Planck Institute for Mathematics in the Sciences, 04103, Leipzig, Germany}
\email{patrick.heslin@mis.mpg.de}
\address{L. Lichtenfelz: Wake Forest University, Winston-Salem, NC., U.S.A.}
\email{lichtel@wfu.edu}
\address{G. Misio{\l}ek: University of Notre Dame, Notre Dame, IN., U.S.A.}
\email{gmisiole@nd.edu}

\begin{document}


\title{Singular sets of Riemannian exponential maps in hydrodynamics}





\begin{abstract}
We prove that the singular sets for the Lagrangian solution maps of the two-dimensional inviscid Euler and generalized surface quasi-geostrophic equations are Gaussian null sets. To achieve this we carry out a spectral analysis of an operator related to the coadjoint representation of the algebra of divergence-free vector fields on the fluid domain. In particular, we establish sharp results on the Schatten-von Neumann class to which this operator belongs. Furthermore, its failure to be compact is directly connected to the absence of Fredholm properties of the corresponding Lagrangian solution maps. We show that, for the three-dimensional Euler equations and the standard surface quasi-geostrophic equation, this failure is in fact essential.
\end{abstract}

\maketitle

\tableofcontents

\section{Introduction}
The Euler equations of ideal hydrodynamics
\begin{gather}\label{eq_euler_intro}
\begin{split}
\partial_t u + \nabla_uu &= -\nabla p \\ 
\div\, u &= 0
\end{split}
\end{gather}
  describe the motion of an incompressible and inviscid fluid. Local existence and uniqueness of solutions was first established in the 1920s by G{\"u}nther and Lichtenstein. Shortly thereafter, Wolibner proved that solutions exist globally in time when the fluid is two-dimensional. Significant improvements were obtained by Yudovich, Kato and others. Most of these developments and the corresponding references can be found for example in the book of Majda \& Bertozzi \cite{majda2002vorticity} and the recent survey of Drivas \& Elgindi \cite{drivas2023singularity}.

In 1966 Arnold \cite{arnold1966sur} introduced a new geometric perspective on the subject. Namely, if the group of volume-preserving diffeomorphisms, viewed as the configuration space of the ideal fluid, is equipped with a right-invariant $L^2$ (kinetic energy) metric, the resulting geodesic equation reduces to the Euler equations above. Ebin \& Marsden \cite{ebin1970groups} proved that the geodesics solve an ordinary differential equation and the associated Riemannian exponential map, which can be viewed as the Lagrangian solution map of \eqref{eq_euler_intro}, is in fact smooth. For a comprehensive overview, see the book of Arnold \& Khesin \cite{arnold2021topological}.

One motivation of Arnold in introducing this geometric picture was to investigate the existence of chaotic attractors in hydrodynamics, hoping that the presence of such a structure may help in understanding the phenomenon of fluid turbulence. Chaotic attractors were known to result from geodesic flows on negatively curved surfaces. However, Arnold found that, while mostly negative, there were directions in which the curvature of the configuration space of the fluid was positive. For this reason, he asked whether conjugate points exist in this setting. Equivalently, one can ask if the set of conjugate vectors of the $L^2$ exponential map, referred to as the singular set, is non-empty. The first non-trivial example was constructed in \cite{misiolek1996conjugate}. Many more examples have followed since then \cite{benn2021conjugate, drivas2022conjugate, lebrigant2024conjugate, preston2006on, shnirelman1994generalized, tauchi2022existence}.

From a fluid dynamics perspective, initial velocity profiles corresponding to conjugate vectors exhibit a type of Lagrangian stability in the sense that the $t=1$ configuration is robust under small perturbations of the initial profile in certain directions. Related to this is the question of uniqueness for the two-point boundary value problem in the Lagrangian framework. In other words, whether one can recover the intermediate flow from two given instances of the fluid's configuration. The answer turns out to be no. Indeed, in any neighbourhood of a conjugate vector, there exist distinct initial velocities whose trajectories arrive at the same configuration at the same time, cf. \cite{lichtenfelz2018normal, misiolek2015the}.

To date, the best known results on the global structure of the singular sets are topological in nature, stemming from the investigation carried out in \cite{ebin2006singularities}. Namely, if the fluid is two-dimensional, the $L^2$ exponential map is a non-linear Fredholm map of index 0, whereas if the fluid is three-dimensional, it is not Fredholm. This yields several consequences for the 2D case. Firstly, the Sard-Smale theorem tells us that the set of conjugate points is of first Baire category - that is, topologically small. Secondly, the set of regular conjugate vectors (which is open and dense among all conjugate vectors) forms a smooth codimension-one submanifold of the tangent space at the identity, cf. \cite{lichtenfelz2018normal}. In contrast, conjugate points in 3D can have infinite order, may accumulate along finite geodesic segments and almost nothing is known about the structure of their corresponding conjugate vectors. However, if the configuration space of a three-dimensional fluid is instead equipped with a Sobolev $H^r$ metric with $r>0$, the corresponding Riemannian exponential map is in fact Fredholm \cite{misiolek2010fredholm}. This suggests that the failure of Fredholmness in the case of 3D ideal fluids $(r=0)$ is, in a certain sense, borderline. Similar studies have been conducted in other settings \cite{li2022vorticity, lichtenfelz2022axisymmetric, shnirelman2005microglobal} including that of the inviscid generalized surface quasi-geostrophic equations \cite{bauer2024geometric}.

The central results of this paper are measure-theoretic in nature. We prove that the singular sets of the exponential maps for the two-dimensional Euler and generalized surface quasi-geostrophic equations (excluding the standard SQG) are Gaussian null - meaning that they have measure zero with respect to any Gaussian measure on the algebra of vector fields. This illustrates the smallness of these singular sets from a measure-theoretic perspective - a notion distinct, even in finite dimensions, from the aforementioned topological smallness. Our proof leverages the analytic dependence of the Lagrangian solution maps on the initial data shown by Shnirelman \cite{shnirelman2012on} and Vu \cite{truong2025thesis}, as well as the existence of an analytic determinant map for Schatten-von Neumann operators obtained by Boyd \& Snigireva \cite{boyd2019on}.

At the heart of establishing Fredholmness in most of the above contexts is the compactness of a certain operator $K$ related to the coadjoint representation of the Lie algebra of vector fields on the underlying domain. In order to establish our results, we first perform a more detailed analysis of the spectrum of this operator. In particular, we show that for the $L^2$ metric this operator in the two-dimensional setting belongs to the Schatten-von Neumann $p$-class for any $p>2$. However, it is not Hilbert-Schmidt and furthermore the failure of compactness of $K$ in the three-dimensional setting is in fact essential. Through explicit examples we demonstrate that the Schatten-$p$ bounds we obtain are sharp, in that there exist initial velocities for which $K$ fails to be Schatten-$p$ for any $p$ outside the established range.

The paper is organized as follows. Section \ref{preliminaries} contains the required preliminaries. Section \ref{the schatten class of K_0} establishes the precise Schatten-von Neumann class of the operator $K$. Section \ref{measures of non-compactness} investigates the critical cases where $K$ fails to be compact. Section \ref{Gaussian nullity of singular sets} contains the proofs of the measure-theoretic results for the singular sets. Appendix \ref{a basis of curl eigenfields} contains the construction of a basis of curl eigenfields on the round three-sphere.

\subsection*{Acknowledgments} 
P. Heslin was supported in part by the National University of Ireland's Dr. {\'E}amon de Valera Postdoctoral Fellowship. L. Lichtenfelz acknowledges support from the A. J. Sterge Faculty Fellowship. G. Misio{\l}ek thanks the Simons Center for Geometry and Physics, where part of this project was carried out, for their hospitality and support. The authors also express their gratitude to Adam Black and Jason Cantarella for helpful conversations and insights for Appendix \ref{a basis of curl eigenfields}.
\newpage

 \section{Preliminaries}\label{preliminaries}
\subsection{Manifold Structure of Diffeomorphism Groups}
Let $M$ be a closed Riemannian manifold and recall the Hodge Laplacian $\Delta = d\delta + \delta d$ where $d$ is the usual exterior derivative and $\delta$ its formal adjoint. Let $s>\frac{\dim M}{2}+1$ (an assumption we will hold throughout the paper) and denote by $H^s(TM)$ the completion of the space of smooth vector fields $\mathfrak{X}(M)$ under the inner product
\begin{equation}\label{full sobolev metric}
    \ip{u,v}_{H^s} = \ip{u,v}_{L^2} + \ip{\Delta^s u, v}_{L^2} \ , \quad u, v \in \mathfrak{X}(M).
\end{equation}
The space $H^s(M,M)$ of self maps of $M$ which are of class $H^s$ in every chart is a smooth Banach manifold. It is modeled on $H^s(TM)$ and continuously embeds into $C^1(M,M)$. The set 
\begin{equation*}
    \D^{s}(M)=\Big\{\eta\in H^{s}(M,M): \ \eta^{-1} \ \text{exists and} \ \eta^{-1}\in H^{s}(M,M)\Big\}
\end{equation*}
of Sobolev diffeomorphisms of $M$ is an open subset of $H^s(M,M)$ and a topological group under composition. The operation of right translation by $\eta \in \D^s(M)$, namely $\xi \mapsto R_\eta \xi = \xi\circ \eta$, is smooth in the $H^s$ topology, but left translation $\xi \mapsto L_\eta \xi = \eta\circ \xi$ is only continuous, as can be seen from their derivatives
\begin{equation*}
    dR_{\eta} v = v\circ\eta \quad \text{and} \quad dL_{\eta} v = D\eta\cdot v, \quad \text{where } v \in T_e\D^s(M).
\end{equation*}
Consequently, operators involving $dL_\eta$ are not generally bounded on $T_{e}\D^{s}(M)$ and will usually be defined on a space of lower-order regularity, $T_{e}\D^{s^{_\prime}}(M)$, the completion of $\mathfrak{X}(M)$ in the $H^{s^{_\prime}}$ norm for $\frac{\dim M}{2} < s^{_\prime}\leq s-1$. For reasons such as this, the space $\D^s(M)$ is strictly speaking not a Lie group in the classical sense. However, it possesses sufficient structure to condone the language we borrow from the theory.

For any $\eta\in\D^{s}\left(M\right)$ consider the Lie group adjoint
\begin{equation*}
    \Ad_{\eta}:\,T_{e}\D^{s^{_\prime}}(M)\rightarrow T_{e}\D^{s^{_\prime}}(M)
\end{equation*}
defined by the usual action of diffeomorphisms on vector fields via pushforward
\begin{equation}\label{DAdjoint}
\Ad_{\eta}v=dL_{\eta}dR_{\eta^{-1}} v= \big(D\eta\cdot v \big)\circ\eta^{-1}=\eta_{*}v.
\end{equation}
Differentiating this expression in $\eta$ gives the Lie algebra adjoint
\begin{equation*}
    \ad_{v}:\,T_{e}\D^{s^{_\prime}}(M)\rightarrow T_{e}\D^{s^{_\prime}}(M)
\end{equation*}
which coincides with the negative of the Lie bracket of vector fields
\begin{equation}\label{AAdjoint}
\ad_{v}w=-[v,w].
\end{equation}

If $\mu$ denotes the volume form on $M$ then the subgroup of volume-preserving $H^s$ Sobolev diffeomorphisms, defined as those preserving $\mu$ under pullback
\begin{equation}
    \Dmu^s(M) = \big\{\eta \in \D^s(M) : \eta^* \mu = \mu\big\},
\end{equation}
is a smooth submanifold of $\D^s(M)$. Its tangent space at the identity consists of divergence-free $H^s$ Sobolev vector fields 
\begin{equation*}
    T_e\Dmu^s(M) = H^s(TM) \cap \mathrm{div}^{-1}(0),
\end{equation*}
while the tangent space at an arbitrary $\eta$ can be described in terms of right-translation
\begin{equation*}
    T_\eta\Dmu^s(M)= \big\{v \circ \eta : v \in T_e\Dmu^s(M)\big\}.
\end{equation*}

Furthermore, $\Dmu^s(M)$ contains the subgroup $\Dex^s(M)$ of exact volume-preserving diffeomorphisms whose Lie algebra is given by
\begin{equation*}
    T_e\Dex^s(M) = \big\{v \in T_e\Dmu^s(M) : \pi_{\mathcal{H}}v = 0\big\}
\end{equation*}
where $\pi_{\mathcal{H}}$ is the projection onto the subspace of harmonic fields in the $L^2$-orthogonal Hodge decomposition, cf. \cite{misiolek2010fredholm}. Note that if $M$ has trivial first cohomology, we have $\Dmu^s(M) = \Dex^s(M)$.

Lastly, if $M$ admits a Killing field $X$ generating a $1$-parameter subgroup of isometries $\{ \Phi_t \}_{t \in \mathbb{R}}$, then the subgroup of axisymmetric diffeomorphisms relative to $X$ is given by
\begin{equation}\label{axisymmetric diffeos}
    \Amu^s(M) = \big\{\eta \in \Dmu^s(M) : \eta \circ \Phi_t = \Phi_t \circ \eta,~\forall t \big\}.
\end{equation}
Its Lie algebra
\begin{equation}\label{eq_axi_Lie_algebra}
    T_e\Amu^s(M) = \big\{v \in T_e\Dmu^s(M) : [v, X] = 0 \big\}
\end{equation}
contains the subspace of swirl-free vector fields, defined as
\begin{equation}\label{swirl free vector fields}
T_e\Swi^s(M) = \big\{v \in T_e\Amu^s(M) : \langle v, X \rangle = 0 \big\}.
\end{equation}
\begin{remark}
Although we write $T_e\Swi^s(M)$ to denote this particular subspace, it need not correspond to the Lie algebra of an actual Lie group. Such a group exists only when the hyperplane distribution induced by $X^{\perp}$ is integrable, cf. \cite{lichtenfelz2022axisymmetric}.
\end{remark}
For technical convenience we will work primarily with $\Dex^s(M)$.

\subsection{Right-invariant Sobolev Metrics}
Given a right-invariant metric $\ip{\cdot, \cdot}$ on a general Lie group $G$ the resulting geodesic flow and its properties can be described in terms of the coadjoint operators
\begin{equation}
    \ip{\Adstar{g} u, v} = \ip{u, \Ad_{g} v} \quad \mathrm{and} \quad \ip{\adstar{w} u, v} = \ip{u, \ad_{w} v}, \quad g \in G ~~\text{and}~ u, v, w \in T_eG
\end{equation}
which combine information from both the Riemannian and Lie structures. The associated geodesic equation when right-translated to the identity $T_eG$ becomes
\begin{equation}\label{euler-arnold}
    \partial_t u = -\adstar{u} u
\end{equation}
and is known as the Euler-Arnold equation. Imposing an initial condition $u(0) = u_0$ and setting $\dot{\gamma} = dR_{\gamma} u$ yields the conservation law
\begin{equation}
    \Adstar{\gamma(t)} u(t) = u_0. 
\end{equation}
Many notable partial differential equations can be written in this fashion, cf. \cite{arnold2021topological, vizman2008geodesicequations}.

As observed by Arnold, the geodesic equation obtained by equipping the space of smooth volume-preserving diffeomorphisms of $M$ with a weak right-invariant $L^2$ metric\footnote{Here we say that a metric is \textit{weak} if it induces a coarser topology than the inherent manifold structure.}, when reduced to the Lie algebra of divergence-free fields, yields the Cauchy problem for the incompressible Euler equations
\begin{equation}\label{euler}
    \begin{split}
        &\partial_t u + \nabla_u u = -\nabla p \\
        &\mathrm{div} \ u = 0 \\
        &u(0) = u_0. 
    \end{split}
\end{equation}
The Riemannian exponential map at the identity is then a natural data-to-solution map for \eqref{euler} expressed in Lagrangian coordinates. A rigorous analysis can be conveniently carried out in the framework of Sobolev $H^s$ spaces where the geodesic equation of the $L^2$ metric becomes an ordinary differential equation on the Banach manifold $\Dmu^s(M)$ and one can apply the standard ODE techniques, cf. \cite{ebin1970groups}. Then, as in classical Riemannian geometry, one has that the $L^2$ Riemannian exponential map
    \begin{equation*}
        \exp_e : U_e \subset T_e\Dmu^s(M) \rightarrow \Dmu^s(M)
    \end{equation*}
is a local diffeomorphism in an open set $U_e$ containing the zero vector. More precisely 
\begin{equation*}
u_0 \mapsto \exp_eu_0 = \gamma(1), 
\end{equation*}
where $\gamma(t)$ satisfies the flow equation 
\begin{equation*}
\dot{\gamma} = u \circ \gamma, \qquad \gamma(0) = e,
\end{equation*}
of the velocity field $u(t)$, solving \eqref{euler}. It is known that if the initial velocity is of higher Sobolev regularity (even smooth), then the geodesic evolves in the space of diffeomorphisms of the same regularity and exists for the same time, cf. \cite{ebin1970groups}.

In this paper we will also consider diffeomorphism groups equipped with right-invariant $H^r$ metrics
\begin{equation}\label{weak full sobolev metric}
    \ip{u,v}_{H^r} = \bigip{(1+\Delta)^ru,v}_{L^2} \quad \text{ for } \ r\leq s-1,
\end{equation}
where $u$ and $v$ are $H^s$ vector fields. If $M$ has trivial first cohomology, or if we restrict to $\Dex^s(M)$, then it is equivalent to work with homogeneous $H^r$ metrics
\begin{equation}\label{homogeneous sobolev metric}
    \ip{u,v}_{\dot{H}^r} = \ip{\Delta^r u, v}_{L^2}.
\end{equation}

The family of the generalized surface quasi-geostrophic equations
\begin{equation}\label{beta SQG}
    \begin{split}
        &\partial_t \theta + u \cdot \nabla \theta = 0 \\
        &u=\symp \Delta^{\frac{\beta}{2}-1}\theta, \qquad 0\leq \beta \leq 1,
    \end{split}
\end{equation}
introduced in \cite{constantin1994formation, cordoba2005evidence} interpolates between the standard SQG equation ($\beta=1$) and the 2D Euler equations ($\beta=0$). These equations are also Euler-Arnold for the metric \eqref{homogeneous sobolev metric} on $\Dex^s(M^2)$ with $r=-\frac{\beta}{2}$ and admit smooth Riemannian exponential maps, cf. \cite{misiolek2023on, truong2025thesis, washabaugh2016the}.

Euler-Arnold equations arising from equipping $\D^s(M)$, $\Dmu^s(M)$ and $\Dex^s(M)$ with the metrics \eqref{weak full sobolev metric} for integer $r\geq0$ include the Euler-$\alpha$ equations ($r=1$) introduced in \cite{holm1998euler}. As shown in \cite{lichtenfelz2022axisymmetric}, this approach also extends to axisymmetric Euler flows. In particular, if $u_0 \in T_e\Amu^s(M)$ then the solution $u(t)$ of the Euler equations belongs to $T_e\Amu^s(M)$ for as long as it exists. Consequently, $\Amu^s(M)$ is a totally geodesic submanifold of $\Dmu^s(M)$ with the induced kinetic energy metric.
In three dimensions, if $u_0$ is in addition swirl-free, then $u(t)$ exists for all time.

The following table illustrates the orders $r$ of Sobolev metrics for which a smooth Lagrangian data-to-solution map is known to exist. Note that all non-negative integers $r=0,1,2\dots$ are included in the first two axes. \\

\begin{tikzpicture}[scale=1.5]
    
    \draw[->] (-1,2.5) -- (7,2.5) node[below] {$r$};
    \node[left] at (-1.3,2.5) {$\mathscr{D}_\mu^s(M^2)$};

    \foreach \x in {-0.5} {
        \draw (\x,2.4) -- (\x,2.6);
        \node[below] at (\x,2.45) {\small $-\frac{1}{2}$};
    }

    \foreach \x in {0.5} {
        \draw (\x,2.4) -- (\x,2.6);
        \node[below] at (\x,2.4) {\small $0$};
    }

    \foreach \x in {2.5} {
        \draw (\x,2.4) -- (\x,2.6);
        \node[below] at (\x,2.4) {\small $1$};
    }

    \foreach \x in {4.5} {
        \draw (\x,2.4) -- (\x,2.6);
        \node[below] at (\x,2.4) {\small $2$};
    }

    \draw[ultra thick,darkgreen] (-0.5,2.5) -- (0.5,2.5);

    \fill[darkgreen] (-0.5,2.5) circle (2pt);
    \fill[darkgreen] (0.5,2.5) circle (2pt);
    \fill[darkgreen] (2.5,2.5) circle (2pt);
    \fill[darkgreen] (4.5,2.5) circle (2pt);

    \draw [decorate,decoration={brace,amplitude=10pt}] (-0.5,2.7) -- (0.5,2.7) node[midway,yshift=0.3cm,above] {\scriptsize Generalized SQG};

    \draw[->] (-0.5,1.95) -- (-0.5,2.05); 
    \node[below] at (-0.5,2) {\scriptsize SQG}; 
    
    \draw[->] (0.5,1.95) -- (0.5,2.05);
    \node[below] at (0.5,2) {\scriptsize 2D Euler};

    \draw[->] (2.5,1.95) -- (2.5,2.05);
    \node[below] at (2.5,2) {\scriptsize 2D Euler-$\alpha$};

    \draw [decorate,decoration={brace,amplitude=10pt}] (4.5,2.7) -- (7,2.7) node[midway,yshift=0.3cm,above] {\scriptsize Higher-order equations};
    
    \draw[->] (-1,1) -- (7,1) node[below] {$r$};
    \node[left] at (-1.3,1) {$\mathscr{D}_\mu^s(M^3)$};

    \foreach \x in {0.5} {
        \draw (\x,0.9) -- (\x,1.1);
        \node[below] at (\x,0.9) {\small $0$};
    }

    \foreach \x in {2.5} {
        \draw (\x,0.9) -- (\x,1.1);
        \node[below] at (\x,0.9) {\small $1$};
    }

    \foreach \x in {4.5} {
        \draw (\x,0.9) -- (\x,1.1);
        \node[below] at (\x,0.9) {\small $2$};
    }

    \draw[->] (0.5,0.45) -- (0.5,0.55);
    \node[below] at (0.5,0.5) {\scriptsize 3D Euler};

    \draw[->] (2.5,0.45) -- (2.5,0.55);
    \node[below] at (2.5,0.5) {\scriptsize 3D Euler-$\alpha$};

    \draw [decorate,decoration={brace,amplitude=10pt}] (4.5,1.2) -- (7,1.2) node[midway,yshift=0.3cm,above] {\scriptsize Higher-order equations};

    \fill[darkgreen] (0.5,1) circle (2pt);
    \fill[darkgreen] (2.5,1) circle (2pt);
    \fill[darkgreen] (4.5,1) circle (2pt);
    
    \draw[->] (-1,-0.5) -- (7,-0.5) node[below] {$r$};
    \node[left] at (-1.3,-0.5) {$\mathscr{A}_\mu^s(M^3)$};

    \foreach \x in {0.5} {
        \draw (\x,-0.6) -- (\x,-0.4);
        \node[below] at (\x,-0.6) {\small $0$};
    }

    \draw[->] (0.5,-1.05) -- (0.5,-0.95);
    \node[below] at (0.5,-1) {\scriptsize Axisymmetric Euler};

    \fill[darkgreen] (0.5,-0.5) circle (2pt);

    \node at (3,-1.75) {\small \parbox{14cm}{\centering
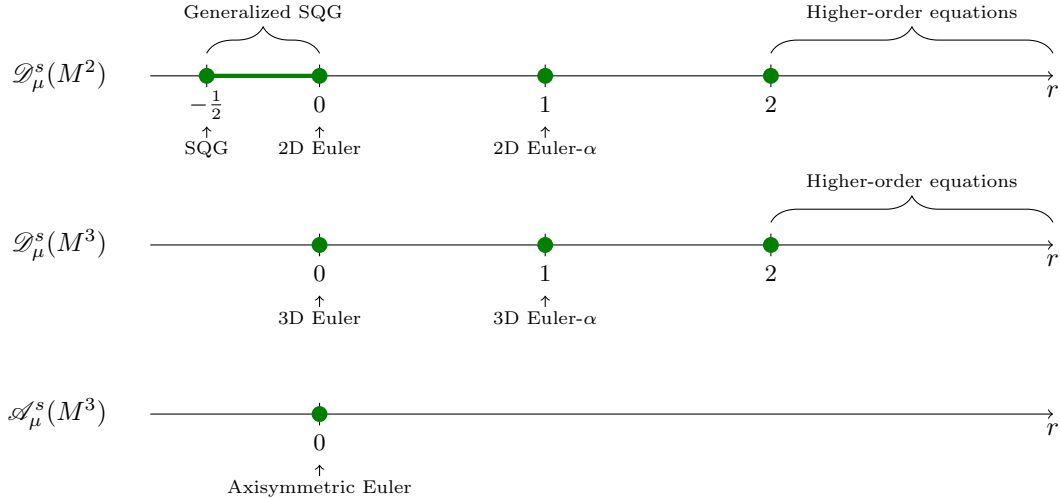
    \textbf{Figure 1: Orders $r$ of the metric \eqref{weak full sobolev metric} for which a smooth Lagrangian solution map is known to exist.}}};
    
\end{tikzpicture}

\subsection{Fredholm Properties of the Exponential Map} Divergence-free vector fields $v$ for which
\begin{equation}\label{dexp}
    d \exp_e(v) : T_e\Dmu^s(M) \rightarrow T_{\exp_e(v)}\Dmu^s(M)
\end{equation}
is not an isomorphism are called monoconjugate if the map \eqref{dexp} fails to be injective and epiconjugate if it fails to be surjective. For exponential maps on general infinite-dimensional manifolds the two types of singularities need not coincide, can be of infinite order and may accumulate along finite geodesic segments. However, the presence of Fredholmness guarantees significantly more structure. As mentioned in the introduction, on $\Dmu^s(M)$ the $L^2$ exponential map is Fredholm of index $0$ if $M$ is two-dimensional but not if $M$ is three-dimensional. For the generalized surface quasi-geostrophic equations \eqref{beta SQG} if $-\frac{1}{2}<r\leq0$ the exponential map is Fredholm of index $0$. If $r=-\frac{1}{2}$ the exponential map fails to be Fredholm. For integer order metrics \eqref{weak full sobolev metric} if $M$ is two-dimensional and $r>-\frac{1}{2}$ or if $M$ is three-dimensional and $r>0$, then the exponential map is again Fredholm of index $0$. These results raise a natural question of whether the failures of Fredholmness for the settings of the standard SQG equation in two-dimensions ($r=-\frac{1}{2}$) and ideal hydrodynamics in three-dimensions ($r=0$) are borderline. We investigate this in Section \ref{measures of non-compactness}.

For the three-dimensional axisymmetric case if $v \in T_e\Amu(M)$ is swirl-free, then $d \exp_e(v) : T_e\Amu^s(M) \rightarrow T_{\exp_e(v)}\Amu^s(M)$ is Fredholm of index 0, cf. \cite{lichtenfelz2022axisymmetric}. By a standard perturbation argument, this remains valid for all axisymmetric velocities with sufficiently small swirl.

The following diagram summarizes the above. Green color indicates that the Fredholm property has been (at least formally) established for this value of $r$. Red dots indicate that the exponential map fails to be Fredholm for this value of $r$. \\

\begin{tikzpicture}[scale=1.5]
    
    \draw[->] (-1,2.5) -- (7,2.5) node[below] {$r$};
    \node[left] at (-1.3,2.5) {$\mathscr{D}_\mu^s(M^2)$};

    \foreach \x in {-0.5} {
        \draw (\x,2.4) -- (\x,2.6);
        \node[below] at (\x,2.45) {\small $-\frac{1}{2}$};
    }

    \foreach \x in {0.5} {
        \draw (\x,2.4) -- (\x,2.6);
        \node[below] at (\x,2.4) {\small $0$};
    }

    \foreach \x in {2.5} {
        \draw (\x,2.4) -- (\x,2.6);
        \node[below] at (\x,2.4) {\small $1$};
    }

    \foreach \x in {4.5} {
        \draw (\x,2.4) -- (\x,2.6);
        \node[below] at (\x,2.4) {\small $2$};
    }

    \draw[ultra thick,darkgreen] (-0.5,2.5) -- (7,2.5);

    \fill[darkred] (-0.5,2.5) circle (2pt);
    
    \draw [decorate,decoration={brace,amplitude=10pt}] (-0.5,2.7) -- (0.5,2.7) node[midway,yshift=0.3cm,above] {\scriptsize Generalized SQG};

    \draw[->] (-0.5,1.95) -- (-0.5,2.05); 
    \node[below] at (-0.5,2) {\scriptsize SQG}; 
    
    \draw[->] (0.5,1.95) -- (0.5,2.05);
    \node[below] at (0.5,2) {\scriptsize 2D Euler};

    \draw[->] (2.5,1.95) -- (2.5,2.05);
    \node[below] at (2.5,2) {\scriptsize 2D Euler-$\alpha$};

    \draw [decorate,decoration={brace,amplitude=10pt}] (4.5,2.7) -- (7,2.7) node[midway,yshift=0.3cm,above] {\scriptsize Higher-order equations};
    
    \draw[->] (-1,1) -- (7,1) node[below] {$r$};
    \node[left] at (-1.3,1) {$\mathscr{D}_\mu^s(M^3)$};

    \foreach \x in {0.5} {
        \draw (\x,0.9) -- (\x,1.1);
        \node[below] at (\x,0.9) {\small $0$};
    }

    \foreach \x in {2.5} {
        \draw (\x,0.9) -- (\x,1.1);
        \node[below] at (\x,0.9) {\small $1$};
    }

    \foreach \x in {4.5} {
        \draw (\x,0.9) -- (\x,1.1);
        \node[below] at (\x,0.9) {\small $2$};
    }

    \draw[->] (0.5,0.45) -- (0.5,0.55);
    \node[below] at (0.5,0.5) {\scriptsize 3D Euler};

    \draw[->] (2.5,0.45) -- (2.5,0.55);
    \node[below] at (2.5,0.5) {\scriptsize 3D Euler-$\alpha$};

    \draw [decorate,decoration={brace,amplitude=10pt}] (4.5,1.2) -- (7,1.2) node[midway,yshift=0.3cm,above] {\scriptsize Higher-order equations};

    \draw[ultra thick,darkgreen] (0.5,1) -- (7,1);

    \fill[darkred] (0.5,1) circle (2pt);
    
    \draw[->] (-1,-0.5) -- (7,-0.5) node[below] {$r$};
    \node[left] at (-1.3,-0.5) {$\mathscr{A}_\mu^s(M^3)$};

    \foreach \x in {0.5} {
        \draw (\x,-0.6) -- (\x,-0.4);
        \node[below] at (\x,-0.6) {\small $0$};
    }

    \draw[->] (0.5,-1.05) -- (0.5,-0.95);
    \node[below] at (0.5,-1) {\scriptsize \shortstack{Axisymmetric Euler \\ (small swirl)}};

    \fill[darkgreen] (0.5,-0.5) circle (2pt);

    \node at (3,-2) {\small \parbox{14cm}{\centering 
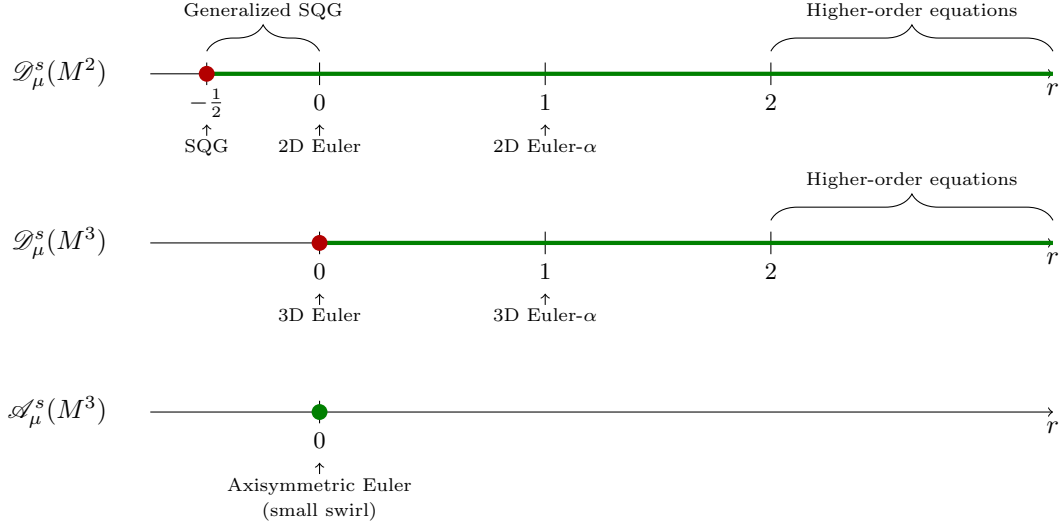
    \textbf{Figure 2: Fredholmness for $H^r$ exponential maps on diffeomorphism groups.}
}};
    
\end{tikzpicture}

The methods of proof for the above results follow similar lines. Fixing smooth\footnote{This allows one to consider the proceeding operators as bounded linear operators on $H^s$.} initial data $u_0 \in T_e\Dex(M)$, consider the corresponding geodesic $t \mapsto \gamma(t) = \exp_e(tu_0)$ in $\Dex(M)$ which, without loss of generality, can be assumed to exist for some time $T>1$. Let $\eta = \gamma(1)$ and consider the operators
\begin{equation*}
    \Lambda(u_0, t): T_e\Dex^s(M) \rightarrow T_e\Dex^s(M), \quad w \mapsto \Adstar{\gamma(t)} \Ad_{\gamma(t)}w
\end{equation*}
and
\begin{equation*}
    K(u_0): T_e\Dex^s(M) \rightarrow T_e\Dex^s(M), \quad w \mapsto \adstar{w} u_0.
\end{equation*}
From these construct
\begin{equation}\label{phi}
    \Phi(u_0, t) : T_e\Dex^s(M) \rightarrow T_e\Dex^s(M), \quad w \mapsto dL_{\gamma(t)}^{-1} d\exp_e(tu_0)tw,
\end{equation}
\begin{equation}\label{omega}
    \Omega(u_0, t) : T_e\Dex^s(M) \rightarrow T_e\Dex^s(M), \quad w \mapsto \int_0^t \Lambda^{-1}(u_0, \tau)w \ d\tau
\end{equation}
and
\begin{equation}\label{gamma}
    \Gamma(u_0, t) : T_e\Dex^s(M) \rightarrow T_e\Dex^s(M), \quad w \mapsto \int_0^t \Lambda^{-1}(u_0, \tau) K(u_0) \Phi(u_0, \tau) w \ d\tau.
\end{equation}
It is not difficult to show that $\Phi$ satisfies the differential equation
\begin{equation}\label{phi ode}
    \frac{d}{dt} \Phi(u_0, t) = \Lambda^{-1}(u_0, t) + \Lambda^{-1}(u_0, t)K(u_0)\Phi(u_0, t),
\end{equation}
which, in light of \eqref{phi}-\eqref{gamma}, when integrated, yields the expression
\begin{equation}\label{dexp decomposition}
    d\exp_e(u_0)w = d_eL_{\eta} \big(\Omega(u_0) w + \Gamma(u_0) w\big).
\end{equation}

The operator $\Omega$ in \eqref{omega} can be shown to be an isomorphism. The key step in proving the Fredholmness result is then showing that, for smooth $u_0$, the operator $K(u_0): T_e\Dex^s(M) \rightarrow T_e\Dex^s(M)$ is compact. A perturbation argument is used to cover the case of $u_0 \in T_e\Dex^s(M)$.

A question left open in \cite{misiolek2010fredholm} concerns the sharpening of these results. More precisely, is $K(u_0)$ Hilbert-Schmidt or trace-class? We show that the answer depends on the value of $r$. To this end we introduce the notation
\begin{equation}\label{K_r(u_0)}
    K_r(u_0): w \mapsto \radstar{w} u_0
\end{equation}
where $\star_r$ refers to the adjoint with the respect to the metric \eqref{weak full sobolev metric} or \eqref{homogeneous sobolev metric}, which should be clear from the context.

If $M$ is two-dimensional, then each $w \in T_e\Dex^s(M)$ can be written as $w = \symp \psi_w$ where $\psi_w$ is a stream function of class $H^{s+1}$. Then, for smooth $u_0 = \symp \psi_{u_0} \in T_e\Dex(M)$ we have
\begin{equation}\label{2D K_0}
\begin{split}
    K_r(u_0)w&= \radstar{\symp \psi_w} \symp \psi_{u_0} \\
    &= \symp \Delta^{-1-r} \big\{\psi_w, \Delta^{1+r}\psi_{u_0}\big\} \\
    &= \Delta^{-\frac{1}{2}-r}\left( \Delta^{\frac{1}{2}+r} \symp \Delta^{-1-r} \big\{\psi_w, \Delta^{1+r}\psi_{u_0}\big\} \right)
\end{split}
\end{equation}
where $\{\psi_1,\psi_2\} = \big<\nabla^{\perp} \psi_1, \nabla\psi_2\big>$ is the Poisson bracket of stream functions and $\ip{\cdot, \cdot}$ denotes the Riemannian metric on $M$.

On the other hand, if $M$ is three-dimensional, then, for $w \in T_e\Dex^s(M)$ and $u_0 \in T_e\Dex(M)$ we have
\begin{equation}\label{3D K_0}
    \begin{split}
    {K_r(u_0)w}&={\curl \Delta^{-1-r}\big[w, \curl \Delta^r u_0\big]} \\
    &= \Delta^{-r}\Big( \curl \Delta^{-1} \big[w, \curl \Delta^r u_0\big]\Big).
    \end{split}
\end{equation}

\begin{remark}
    Of course \eqref{2D K_0} and \eqref{3D K_0} remain valid for $u_0 \in T_e\Dex^s(M)$ by considering $K_0$ as a bounded linear operator on a space of lower regularity.
\end{remark}

It may be useful at this point to take note of the difference in powers of the Laplacian in the two and three-dimensional settings. This will play a role later.

\subsection{Schatten-von Neumann Classes}
Let $\mathcal{H}$ be a real separable Hilbert space and $T$ a compact operator on $\mathcal{H}$. The singular values $\{s_n(T)\}_{n \in \mathbb{N}}$ of $T$ are defined as the square roots of the eigenvalues of $T T^\star$ and satisfy $\displaystyle \lim_{n \to \infty} s_n(T) = 0$, where $\star$ denotes the adjoint with respect to the inner product $\ip{\cdot, \cdot}_{\mathcal{H}}$.
\begin{remark}\label{eigenvalue to singular value}
    If $T^\star = \pm T$ then the singular values of $T$ are equal to the absolute values of the eigenvalues of the extension of $T$ to the complexification $\C \otimes \mathcal{H}$. This fact will be useful for several examples we consider in Sections \ref{the schatten class of K_0} and \ref{measures of non-compactness}.
\end{remark}
For $1 \leq p < \infty$, the Schatten-von Neumann $p$-class (or Schatten $p$-class for brevity)
\begin{equation}
\mathcal{S}_p(\mathcal{H}) = \big\{ T : \mathcal{H} \rightarrow \mathcal{H} \, :  \, T \text{ is compact and } \{ s_n(T) \}_{n \in \mathbb{N}} \in \ell^p \big\}
\end{equation}
is a Banach space with the norm
\begin{equation}
\| T \|_{\mathcal{S}_p} = \| s_n(T) \|_p = \left( \sum\limits_{n = 1}^{\infty} s_n(T)^p \right)^{\frac{1}{p}}.
\end{equation}
Special cases include the well-known Hilbert-Schmidt ($p=2$) and trace-class operators ($p=1$). A comprehensive treatment can be found for example in \cite{gohberg1969introduction, schatten1960norm}.

\begin{remark}
    The Schatten $p$-norm can equivalently be computed using an orthonormal basis $\{v_j\}_{j \in \mathbb{N}}$ for $\mathcal{H}$. In particular
    \begin{equation}\label{brute schatten norm}
        \norm{T}_{\mathcal{S}^p}^p = \sum_{j = 1}^\infty \norm{T v_j}_{\mathcal{H}}^p.
    \end{equation}
\end{remark}

As with standard $\ell^p$ spaces, there are continuous embeddings $\mathcal{S}_p(\mathcal{H}) \subseteq \mathcal{S}_q(\mathcal{H})$ whenever $p \leq q$. Note that a compact operator does not necessarily belong to $\mathcal{S}_p(\mathcal{H})$ for any $p$. Each class forms a two-sided ideal within the space $L(\mathcal{H})$ of bounded linear operators on $\mathcal{H}$. In particular, given $T \in \mathcal{S}_p(\mathcal{H})$, and $A, B \in L(\mathcal{H})$, we have the estimate
\begin{equation}\label{eq_two_sided_ideal}
\|ATB\|_{\mathcal{S}_p} \leq \|A\|_{\mathrm{op}}\|T\|_{\mathcal{S}_p}\|B\|_{\mathrm{op}},
\end{equation}
where $\|\cdot \|_{\mathrm{op}}$ is the operator norm. Note that $\mathcal{S}_p(\mathcal{H})$ is not a closed subspace of $L(\mathcal{H})$ for the same reason that $\ell^p$ is not closed in $\ell^{\infty}$. As a consequence, the integral of a curve of Schatten-$p$ operators might fail to be Schatten-$p$, since the passage from Riemann sums to the integral involves taking a limit. Nevertheless, we have the following lemma.

\begin{lemma}\label{schatten_lemma}
Fix $p \in [1, \infty)$ and $T \in \mathcal{S}_p(\mathcal{H})$. Let $A, B : [0, 1] \rightarrow L(\mathcal{H})$ be two continuous one-parameter families of bounded linear operators satisfying uniform estimates $\|A(t)\|_{\mathrm{op}} \leq C_1$ and $\|B(t)\|_{\mathrm{op}} \leq C_2$ for all $t \in [0, 1]$. Then, the integral
\begin{equation}
\mathcal{I} = \int\limits_0^1 A(t)TB(t)\,dt
\end{equation}
defines a Schatten-$p$ class operator.
\end{lemma}

\begin{proof}
This follows by a direct estimate using the uniform bounds and the two-sided ideal property \eqref{eq_two_sided_ideal} of $\mathcal{S}_p(\mathcal{H})$. In particular, one has
\begin{align*}
    \|\mathcal{I}\|_{\mathcal{S}_p} &\leq \int\limits_0^1 \norm{A(t)TB(t)}_{\mathcal{S}_p}\,dt \\
    &\leq \int\limits_0^1 \norm{A(t)}_{\mathrm{op}}\norm{T}_{\mathcal{S}_p}\norm{B(t)}_{\mathrm{op}}\,dt \\
    &\leq C_1C_2\norm{T}_{\mathcal{S}_p}.
\end{align*}
\end{proof}

We will make use of a result due to Boyd \& Snigivera in \cite{boyd2019on}. Namely, for each $p \in [1, \infty)$ there exists a determinant function $\mathrm{det}_p : \big(\mathcal{S}_p(\mathcal{H}), \norm{\cdot}_{\mathcal{S}_p}\big) \rightarrow \mathbb{R}$ with the property that, for any $T \in \mathcal{S}_p(\mathcal{H})$, we have\footnote{The right-hand side of \eqref{eq_det_p} is occasionally written in the literature by an abuse of notation as $\mathrm{det}_p(\mathrm{Id} + T) \neq 0$.}
\begin{equation}\label{eq_det_p}
\mathrm{Id} + T ~\text{is invertible} ~\Leftrightarrow~\mathrm{det}_p(T) \neq 0.
\end{equation}
It will be of importance that this map is in fact analytic.

\subsection{Gaussian Null Sets}\label{gaussian null sets} Let $\mathcal{X}$ be a real separable Banach space and $\mathcal{B}(\mathcal{X})$ its collection of Borel sets. A non-degenerate Gaussian measure $\nu$ on $\mathcal{B}(\mathcal{X})$ is a probability measure such that for all non-zero elements $f$ in the dual $\mathcal{X}^*$ and $B \in \mathcal{B}(\R)$ the induced measure on the real line $\nu \circ f^{-1}$ has the form
\begin{equation*}
    \nu \circ f^{-1}(B) = \frac{1}{\sqrt{2 \pi \beta}} \int_B e^{- \frac{(x-\alpha)^2}{2\beta}} dx
\end{equation*}
for some $\alpha \in \R$ and $\beta>0$ with $dx$ denoting the standard Lebesgue measure.

A Borel subset of $\mathcal{X}$ is called Gaussian null if it has measure zero with respect to any non-degenerate Gaussian measure on $\mathcal{X}$. This notion was introduced by Phelps \cite{phelps1978gaussian} as a generalization of Lebesgue null sets in order to prove an infinite-dimensional analogue of Rademacher's theorem. Notably, every Borel subset of a Gaussian null set is itself Gaussian null. Moreover, a result due to Bogachev and Malofeev \cite{bogachev2014onthe} implies that the zero set of any analytic function on $\mathcal{X}$ is Gaussian null. This will play a key role  later.





\section{The Schatten-von Neumann Class of \texorpdfstring{$K_r(u_0)$}{K0}}\label{the schatten class of K_0}
We now examine to which Schatten $p$-class the operator $K_r(u_0)$ defined in \eqref{K_r(u_0)} belongs. It turns out that the answer depends on the dimension of the underlying manifold and the value of $r$. For this reason we will consider the two and three-dimensional cases separately. Our argument draws on techniques from complex analysis. The results could also be established using real interpolation methods, cf. \cite{cwikel1977weak}.

For $s>\frac{\dim M}{2}+1$, consider the group $\Dex^s$ of exact volume-preserving Sobolev diffeomorphisms of a closed oriented Riemannian manifold $M$. As mentioned in the preliminaries, due to loss of derivatives, one can either assume that $u_0$ is smooth, in which case $K_r(u_0)$ is a bounded linear operator on $T_e\Dex^s$, or that $u_0$ is of class $H^s$, in which case $K_r(u_0)$ is bounded in a coarser topology. For reasons which will be clear in Section \ref{Gaussian nullity of singular sets}, we take the latter approach.


\subsection{The two-dimensional case} Let $M$ be a closed surface. The main result of this subsection is the following.

\begin{theorem}\label{main_thm_2d}
    If $r>-\frac{1}{2}$ and $1 < s^{_\prime} \leq s - 2 - 2r$, then for any $u_0 \in T_e\Dex^s(M)$ the operator $K_r(u_0)$ on $T_e\Dex^{s^{_\prime}}(M)$ given in \eqref{2D K_0} is of Schatten-von Neumann $p$-class for all $p>\frac{2}{1+2r}$.
\end{theorem}

\begin{proof}
    Recall from \eqref{2D K_0} that, for $w =\symp \psi_w \in T_e\Dex^{s^{_\prime}}$, we have
    \begin{equation*}
        K_r(u_0)w = \Delta^{-\frac{1}{2}-r}\left( \Delta^{\frac{1}{2}+r} \symp \Delta^{-1-r} \{\psi_w, \Delta^{1+r}\psi_{u_0}\} \right). 
    \end{equation*}
    For fixed $u_0 \in T_e\Dex^s$ the linear operator in parentheses above, being a composition of a pseudodifferential operator of order $0$ and a multiplication operator with $H^{s-2-2r}$ coefficients, is bounded on $T_e\Dex^{s^{_\prime}}$. In particular we have
    \begin{align*}
        \norm{\Delta^{\frac{1}{2}+r} \symp \Delta^{-1-r} \{\psi_w, \Delta^{1+r}\psi_{u_0}\}}_{H^{s^{_\prime}}} &\lesssim \norm{\big<w, \nabla \Delta^{1+r}\psi_{u_0}\big>}_{H^{s^{_\prime}}} \\
        &\lesssim \norm{\nabla \Delta^{1+r}\psi_{u_0}}_{H^{s^{_\prime}}} \norm{w}_{H^{s^{_\prime}}} \\
        &\lesssim \norm{\nabla \Delta^{1+r}\psi_{u_0}}_{H^{s-2r-2}} \norm{w}_{H^{s^{_\prime}}} \\
        &\lesssim \norm{w}_{H^{s^{_\prime}}}.
    \end{align*}    
    Hence, as each Schatten $p$-class is an ideal in the space of bounded linear operators, a finite Schatten $p$-norm for the inverse of the Laplacian on $T_e\Dex^{s^{_\prime}}$ will imply the same for $K_r(u_0)$.
    
    The Schatten $p$-class of the fractional Laplacian can be determined by the region where the corresponding zeta function 
    is holomorphic, for which there is a large literature. To make the relationship explicit, consider an $H^{s^{_\prime}}\!$-orthonormal basis of eigenfields $\{w_{j}\}_{j \in J}$ of the Hodge Laplacian  with corresponding eigenvalues $\{\lambda_{j}\}_{j \in J}$ for some index set $J$. The zeta function $\zeta(z)$ is then the meromorphic extension of $z \mapsto \sum_{j \in J} \abs{\lambda_j}^{-z},$ which converges for $\mathrm{Re}(z)>\frac{\dim M}{2}$, cf. \cite[Theorem~5.2]{rosenberg1997the}. Furthermore, as $\Delta$ is formally self-adjoint, for any $\sigma>0$, by orthonormality of the basis $\{ w_j \}_{j \in J}$, we have
    \begin{align*}
        \norm{\Delta^{-\sigma}}_{\mathcal{S}_p}^p &= \sum_{j \in J} \norm{\Delta^{-\sigma}w_j}_{H^{s^{_\prime}}}^p = \sum_{j \in J} \abs{\lambda_j}^{-\sigma p} = \zeta(\sigma p)
    \end{align*}
    for $\sigma p>\frac{\dim M}{2} = 1$. The theorem follows by setting $\sigma=\frac{1}{2}+r$.
\end{proof}

Next, we show by explicit example that the above result is sharp.

\begin{theorem}\label{2D Example}
Under the assumptions of Theorem \ref{main_thm_2d} there exists $u_0 \in T_e\Dmu^s(\mathbb{S}^2)$ such that for all $p \leq \frac{2}{1+2r}$ the operator $K_r(u_0)$ is not of Schatten-von Neumann $p$-class on $T_e\Dmu^{s^{_\prime}}(\mathbb{S}^2)$.
\end{theorem}

\begin{proof}
In standard spherical coordinates $(\theta, \phi)$ the round metric is $ds^2 = \sin^2\phi \ d\theta^2 + d\phi^2$. For $f, g$ smooth functions on $\mathbb{S}^2$, we have the following formulas for the skew-gradient
\begin{equation*}
    \symp f = \frac{1}{\sin{\phi}}\Big(\partial_{\phi}f\,  \partial_\theta - \partial_{\theta}f\,
    \partial_\phi\Big),
\end{equation*}
the Poisson bracket
\begin{equation*}
    \{f,g\} = \frac{1}{\sin{\phi}}\left(\partial_{\theta}f\partial_{\phi}g -  \partial_{\phi}f\partial_{\theta} g\right)
\end{equation*}
and the Hodge Laplacian
\begin{equation*}
    \Delta f = -\frac{\partial_{\theta}^2f}{\sin^2 \phi} - \partial_{\phi}^2f - \cot \phi \ \partial_{\phi}f.
\end{equation*}
Consider now the complex spherical harmonics
\begin{equation*}
    Y^m_{\ell}(\theta, \phi) = c^m_{\ell} P^m_{\ell}(\cos \phi) e^{im\theta}, \qquad {\ell} \in \Z_{\geq0}, \ \abs{m} \leq {\ell}
\end{equation*}
where
\[P^m_{\ell}(x)= \frac{(-1)^m}{2^{\ell} {\ell}!}(1-x^2)^{\frac{m}{2}}\frac{d^{m+{\ell}}}{dx^{m+{\ell}}}\left((x^2-1)^{\ell}\right)\]
are the associated Legendre polynomials and $c^m_{\ell}$ are normalizing constants such that $\norm{\symp Y^m_{\ell}}_{\dot{H}^{s^{_\prime}}}=1$. These functions satisfy the eigenvalue equations
\begin{equation}
\label{eq_eigen1}
    \Delta Y^m_{\ell} = {\ell}({\ell}+1)Y^m_{\ell} \qquad \text{and} \qquad \partial_\theta Y^m_{\ell} = imY^m_{\ell}.
\end{equation}
Let $u_0$ be the equatorial rotation vector field $\partial_\theta = \symp Y_1^0$. From \eqref{2D K_0}, we have
\begin{gather}
\begin{split}
\label{eq_eigen2}
    K_r(\partial_\theta)w &= \symp \Delta^{-1-r} \big\{\psi_w, \Delta^{1+r}Y_1^0\big\} \\
    &= 2^{1+r} \symp \Delta^{-1-r}\partial_\theta \psi_w.
\end{split}
\end{gather}
\begin{remark}
When extended to $\C \otimes T_e\Dmu^{s^{_\prime}}$, the complexification of $T_e\Dmu^{s^{_\prime}}$, the operator $K_r(\partial_\theta)$ has $\{\symp Y_{\ell}^m\}$ with $\ell \geq 1$ and $|m| \leq \ell$ as an orthonormal eigenbasis. It is not difficult to check that it is skew self-adjoint with respect to the $H^{s^{_\prime}}$ inner product. Hence its singular values coincide with the absolute values of the purely imaginary eigenvalues of its extension, cf. Remark \ref{eigenvalue to singular value}.
\end{remark}
Using \eqref{eq_eigen1} and \eqref{eq_eigen2}, we obtain the eigenvalues
\begin{equation}\label{rotation eigenvalues}
    K_r(\partial_\theta) \symp Y_\ell^m = im\left(\frac{2}{\ell(\ell + 1)}\right)^{1+r} \symp Y_\ell^m, \quad \text{for} \ \ell \geq 1 \ \text{and} \ \abs{m} \leq \ell.
\end{equation}
Hence, the Schatten-$p$ norm of $K_r(\partial_\theta)$ is
\begin{align*}
    \norm{K_r(\partial_\theta)}_{\mathcal{S}_p}^p &= \sum_{{\ell}=1}^\infty \sum_{m=-{\ell}}^{{\ell}} \abs{m}^p \left(\frac{2}{{\ell}({\ell}+1)}\right)^{p+rp} \\
    &\simeq \sum_{{\ell}=1}^\infty {\ell}^{p+1} \left(\frac{2}{{\ell}({\ell}+1)}\right)^{p+rp}\\
    &\simeq \sum_{{\ell}=1}^\infty \frac{1}{{\ell}^{p+2rp-1}}
\end{align*}
which diverges for $p \leq \frac{2}{1+2r}$.
\end{proof}

\begin{remark}
    We see from the above that for $2$D  hydrodynamics ($r = 0$) the operator $K_0(u_0)$ may fail to be Hilbert-Schmidt. However, the failure is borderline in the sense that it is Schatten $p$-class for any $p>2$. For the standard surface quasi-geostrophic equation ($r = -1/2$), the operator $K_{\scriptscriptstyle{-\frac{1}{2}}}(u_0)$ may even fail to be compact cf. \cite{bauer2024geometric, washabaugh2016the}.
\end{remark}

\subsection{The three-dimensional case} Let $M$ be a closed three-dimensional Riemannian manifold. The main result of this subsection is the following.

\begin{theorem}\label{main_thm_3d}
    If $r>0$ and $\frac{5}{2} < s^{_\prime} \leq s - 2 - 2r$, then for any $u_0 \in T_e\Dex^s(M)$ the operator $K_r(u_0)$ on $T_e\Dex^{s^{_\prime}}(M)$ given in \eqref{3D K_0} is of Schatten-von Neumann $p$-class for all $p>\frac{3}{2r}$.
\end{theorem}

\begin{proof}
Let $w \in T_e\Dex^{s^{_\prime}}$. Recall that we have
\begin{equation*}
    {K_r(u_0)w} = \Delta^{-r}\Big( \curl \Delta^{-1} \big[w, \curl \Delta^r u_0\big]\Big).
\end{equation*}
For $u_0 \in T_e\Dex^s$, the operators in parentheses above are bounded on $T_e\Dex^{s^{_\prime}}$, since they are of zeroth order in $w$. Therefore, as in Theorem \ref{main_thm_2d}, the Schatten $p$-class of $K_r(u_0)$ on $T_e\Dex^{s^{_\prime}}$ is determined by the Schatten $p$-class of $\Delta^{-r}$ on $T_e\Dex^{s^{_\prime}}$. Using the same analysis of the corresponding zeta function, we find that $K_r(u_0)$ is Schatten-$p$ if $rp > \frac{3}{2}$, as claimed.
\end{proof}

\begin{remark}
Despite similarities in the statements, it is important to note that Theorem \ref{main_thm_3d} does not reduce to our previous Theorem \ref{main_thm_2d} when $\dim M = 2$, cf. equations \eqref{2D K_0}, \eqref{3D K_0} and the subsequent remark. Note also the more restrictive condition on $s^{_{\prime}}$. This is a consequence of the difference in dimensions and the fact that the operator $w \mapsto [w, \curl \Delta^r u_0]$ is of order one, whereas $w\mapsto \big\{\psi_w, \Delta^{1+r}\psi_{u_0}\big\}$ is order zero. 
\end{remark}

As in the two-dimensional case, we show by explicit example that the above result is sharp.

\begin{theorem}\label{3D Example}
Under the assumptions of Theorem \ref{main_thm_3d}, there exists $u_0 \in T_e\Dmu^s(\mathbb{S}^3)$ such that for all $p \leq \frac{3}{2r}$ the operator $K_r(u_0)$ is not of Schatten-von Neumann $p$-class on $T_e\Dmu^{s^{_\prime}}(\mathbb{S}^3)$.
\end{theorem} 

\begin{proof} 
Consider the basis for right-invariant vector fields given by
\begin{align*}
e_1 &= -y \partial_x + x \partial_y - w \partial_z + z \partial_w, \\
e_2 &= -z\partial_x + w\partial_y + x\partial_z -y \partial_w, \\
e_3 &= -w\partial_x -z\partial_y + y\partial_z + x\partial_w.
\end{align*}
We compute the spectrum of $K_r(u_0)$ for the initial data $u_0 = e_1$. Since this is a Hopf field, and therefore a curl eigenfield with $\curl e_1 = 2 e_1$, we can see from \eqref{3D K_0} that
\begin{equation}\label{Hopf field K}
    K_r(e_1)w = 2^{2r+1} \curl \Delta^{-1-r}[w, e_1].
\end{equation}
Notice that the kernel of $K_r(e_1)$ consists precisely of those vector fields which are axisymmetric with respect to the Hopf field $e_1$, cf.   \eqref{eq_axi_Lie_algebra}.

As in the two-dimensional case, $K_r(e_1)$ is skew self-adjoint with respect to the $H^{s^{_\prime}}$ metric. Hence, we have again that its singular values coincide with the absolute values of the purely imaginary eigenvalues of its complexification, cf. Remark \ref{eigenvalue to singular value}. From \eqref{Hopf field K} we see that a simultaneous eigenbasis for the Laplacian and the Lie derivative on vector fields $\mathcal{L}_{e_1}: w \mapsto -[w,e_1]$ will give us an eigenbasis for $K_r(e_1)$. In order to count multiplicities, we construct this basis explicitly. On account of the technical nature of the construction, it is presented in Appendix \ref{a basis of curl eigenfields}. We describe the relevant properties of the basis here.

\begin{restatable}{proposition}{basisprop}\label{the basis proposition} Let $e_1 = -y \partial_x + x \partial_y - w \partial_z + z \partial_w$ be the Hopf field on $\mathbb{S}^3$. There exists a Schauder basis of $\C \otimes T_e\Dmu(\mathbb{S}^3)$ consisting of two families of divergence-free vector fields $\mathcal{E}_k^m$ for $k\geq0$ and $0 \leq m \leq k$ and $\mathcal{F}_{k}^m$ for $k \geq 2$ and $1 \leq m \leq k-1$ such that
    \begin{itemize}
    \item The cardinality $\# \mathcal{E}_0^0 = 3$ and, if $v \in \mathcal{E}_0^0$, then $\curl v = 2v$. \\
        \item For $k\geq 1$ the cardinality $\# \mathcal{E}_k^m =
            \begin{cases}
                2(k+1) \quad \text{if } \ m=0, \\
                k+1 \quad \text{if } \ 1\leq m \leq k-1, \\
                2(k+1) \quad \text{if } \ m=k.
            \end{cases}$ \\
        \item For $k\geq1$ if $v \in \mathcal{E}_k^m$ then $\curl v = (k+2)v$ and $[v,e_1] = -i(2m-k)v$. \\
        \item For $k \geq 2$ the cardinality $\# \mathcal{F}_{k}^m = k+1$. \\
        \item For $k \geq 2$ if $v \in \mathcal{F}_{k}^m$ then $\curl v = -kv$ and $[v,e_1] = -i(2m-k)v$.
    \end{itemize}
\end{restatable}

From here we compute that, for $k \geq 1$ and $v \in \mathcal{E}_k^m$ we have
\begin{equation}\label{Hopf field eigenvalues on E_k^m}
    K_r(e_1)v = -i (2m-k)\left(\frac{2}{k+2}\right)^{2r+1}v \quad \text{for} \ 0 \leq m \leq k. 
\end{equation}
while for $k \geq 2$ and $w \in \mathcal{F}_k^m$ we have
\begin{equation}\label{Hopf field eigenvalues on F_k^m}
    K_r(e_1)w = i (2m-k)\left(\frac{2}{k}\right)^{2r+1}w \quad \text{for} \ 1 \leq m \leq k-1. 
\end{equation}
\begin{remark}
    It follows immediately from examining the asymptotics of the eigenvalues above that the operator $K_r(e_1): T_e\Dex^{s^{_\prime}} \rightarrow T_e\Dex^{s^{_\prime}}$ is compact if and only if $r>0$. This is reminiscent of the threshold for the Fredholm property, cf. the introduction.
\end{remark}
We now compute the Schatten-$p$ norm. Ignoring the inconsequential contributions of $v \in \mathcal{E}_0^0$, we have 
\begin{align*}
    \norm{K_r(e_1)}_{\mathcal{S}^p}^p &\simeq \sum_{k=1}^\infty \sum_{m=0}^k \sum_{v \in \mathcal{E}_k^m} \left|-i (2m-k)\left(\frac{2}{k+2}\right)^{2r+1}\right|^p + \sum_{k=2}^\infty \sum_{m=1}^{k-1} \sum_{w \in \mathcal{F}_k^m} \left|-i (2m-k)\left(\frac{2}{k}\right)^{2r+1}\right|^p \\
    &\simeq \sum_{k=1}^\infty \sum_{m=0}^k (k+1) \abs{2m-k}^p\left(\frac{2}{k+2}\right)^{2rp+p} + \sum_{k=2}^\infty \sum_{m=1}^{k-1} (k+1)\abs{2m-k}^p\left(\frac{2}{k}\right)^{2rp+p} \\
    &\simeq \sum_{k=1}^\infty (k+1)k^{p+1}\left(\frac{2}{k+2}\right)^{2rp+p} \\
    &\simeq \sum_{k=1}^\infty \left(\frac{1}{k}\right)^{2rp - 2}.
\end{align*}
Hence, $K_r(e_1)$ fails to be of Schatten $p$-class when $2rp-2\leq1$, that is, if $p\leq\frac{3}{2r}$.
\end{proof}

The above results are summarized graphically below. Green color indicates that, for the given value of $r$ and for any initial data $u_0$, the operator $K_r(u_0): T_e\Dex^{s^{_\prime}} \rightarrow T_e\Dex^{s^{_\prime}}$ belongs to the given Schatten $p$-class. Red color indicates that, for the given value of $r$, there exists an initial data $u_0$ such that $K_r(u_0)$ does not belong to this class.
\medskip

\begin{center}
{\begin{minipage}{0.45\textwidth}
    \centering
    \begin{tikzpicture}
        \begin{axis}[
            axis lines=middle,
            axis on top,
            xlabel={$r$},
            xlabel style={at={(axis description cs:1,-0.075)}},
            ylabel={$p$},
            xmin=-0.5, xmax=5,
            ymin=0, ymax=5.5,
            samples=100,
            domain=-0.5:5,
            axis equal image,
            xtick={0, 1, 3},
            extra x ticks={-0.5},
            extra x tick labels={$-\tfrac{1}{2}$},
            extra x tick style={grid=major}
        ]
        
        \addplot [name path=f, red, thick] {2/(1+2*x)};
        \addplot [name path=axis, draw=none] coordinates {(-0.5,0) (5,0)};
        \addplot [red!20] fill between[of=f and axis];

        \addplot [name path=top, draw=none] coordinates {(-0.5,5.5) (5,5.5)};
        \addplot [green!20] fill between[of=f and top];

        \addplot [dashed, thick] coordinates {(-0.5,0) (-0.5,5.5)};

        \node at (axis cs:3,1) [anchor=north west, fill=green!20, text=red, font=\small] {$p = \frac{2}{1+2r}$};

        \end{axis}
    \end{tikzpicture}
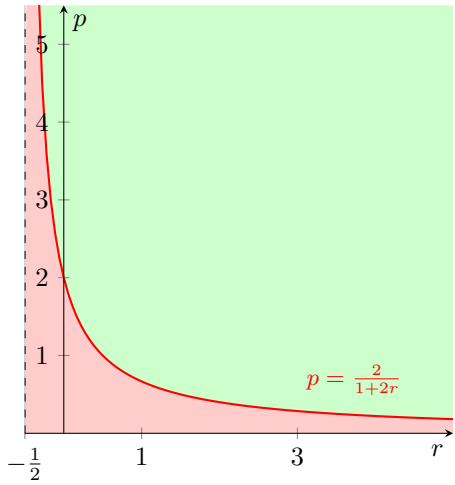
    \raisebox{0.15ex}[0pt][0pt]{\textbf{Two-dimensional setting}}  
\end{minipage}
\hfill
\begin{minipage}{0.45\textwidth}
    \centering
    \begin{tikzpicture}
        \begin{axis}[
            axis lines=middle,
            axis on top,
            xlabel={$r$},
            xlabel style={at={(axis description cs:1,-0.075)}},
            ylabel={$p$},
            ylabel style={at={(axis description cs:-0.01,1)}},
            xmin=0, xmax=5.5,
            ymin=0, ymax=5.5,
            samples=100,
            domain=0:5.5,
            axis equal image,
            xtick={0, 1, 3},
            extra x ticks={0},
            extra x tick labels={0},
            extra x tick style={grid=major}
        ]
        
        \addplot [name path=f, red, thick] {3/(2*x)};
        \addplot [name path=axis, draw=none] coordinates {(0,0) (5.5,0)};
        \addplot [red!20] fill between[of=f and axis];

        \addplot [name path=top, draw=none] coordinates {(0,5.5) (5.5,5.5)};
        \addplot [green!20] fill between[of=f and top];

        \addplot [dashed, thick] coordinates {(0,0) (0,5.5)};

        \node at (axis cs:4,1) [anchor=north west, fill=green!20, text=red, font=\small] {$p = \frac{3}{2r}$};

        \end{axis}
    \end{tikzpicture}
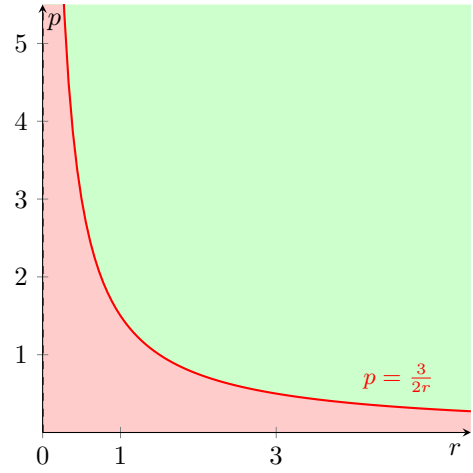
    \raisebox{-0.5ex}[0pt][0pt]{\textbf{Three-dimensional setting}}  
\end{minipage}
}
\end{center}

In fact, Theorem \ref{main_thm_3d} generalizes to higher dimensional manifolds.

\begin{theorem}\label{main_thm_nd}
    Let $\dim M \geq 4$. If $r>0$ and $\frac{\dim M}{2} + 1 < s^{_\prime} \leq s - 2 - 2r$, then for any $u_0 \in T_e\Dex^s(M)$ the operator $K_r(u_0)$ on $T_e\Dex^{s^{_\prime}}(M)$ given in \eqref{3D K_0} is of Schatten-von Neumann $p$-class for all $p>\frac{\dim M}{2r}$.
\end{theorem}

\begin{proof} Notice that, in the language of differential forms
\begin{equation*}
    {K_r(u_0)w^{\flat}} = 
    \delta \Delta^{-1-r}d\iota_{w}\big(d\Delta^ru_0^{\flat}\big)= \Delta^{-r} \Big(\delta \Delta^{-1}d\iota_{w}\big(d\Delta^ru_0^{\flat}\big)\Big).
\end{equation*}
and hence the proof is identical to that of Theorem \ref{main_thm_3d} on account of the common structure. 
\end{proof}

It is reasonable to expect that this is also a sharp result in the same sense as Theorems \ref{2D Example} and \ref{3D Example}.

\subsection{The axisymmetric setting} Recall from the preliminaries that, on a general Riemannian manifold $M$ equipped with a Killing field $X$, we have a notion of axisymmetric vector fields, namely those that are divergence-free and Lie commute with $X$, and the corresponding diffeomorphisms, cf. \eqref{axisymmetric diffeos}-\eqref{eq_axi_Lie_algebra}. This subsection is devoted to the spectral analysis of the operator $K_0(u_0)$ when restricted to axisymmetric vector fields, which is possible due to the fact that $K_0(u_0)$ maps $T_e\Amu^{s^{_\prime}}$ to itself, cf. \cite{lichtenfelz2022axisymmetric}.

We concentrate on the case $M = \mathbb{S}^3$ for concreteness, and take our Killing vector field $X$ to be the Hopf field $e_1$. The axisymmetric condition is then given by $[v,e_1] = 0$. Hence, $\C \otimes T_e\Amu^{s^{_\prime}}$ is spanned by those $\mathcal{E}_k^m$ and $\mathcal{F}_k^m$ from Proposition \ref{the basis proposition} whose elements Lie-commute with $e_1$. If $k$ is odd, no such sets exists. If $k$ is even, say $k=2\ell$, then there are precisely two sets whose elements possess this property: $\mathcal{E}_{2\ell}^\ell$ and $\mathcal{F}_{2\ell}^\ell$. Including the only relevant element from the basis for $k=0$, we have a basis for $\C \otimes T_e\Amu^{s^{_\prime}}$ given by $\{e_1\} \cup \{\mathcal{E}_{2\ell}^\ell\}_{\ell\geq1} \cup \{\mathcal{F}_{2\ell}^\ell\}_{\ell \geq 1}$.

Next, consider any $u_0 \in T_e\Swi^{s}$ swirl-free initial data, i.e. $\ip{u_0,e_1}$ $=$ $0$, cf. \eqref{swirl free vector fields}. From this assumption we have that $\curl u_0 = \phi e_1$ for some function $\phi$ of class $H^{s-1}$. Hence, for any $v$ in our basis we have
\begin{equation*}
K_0(u_0) v = \curl^{-1}[v, \curl u_0] = \curl^{-1} d\phi(v).
\end{equation*}
From this we can estimate for $v \in \mathcal{E}_{2\ell}^\ell \cup \mathcal{F_{2\ell}^\ell}$ 
\begin{equation*}
\|K_0(u_0) v\|_{H^{s^{_\prime}}} \lesssim \| d\phi(v) \|_{H^{{s^{_\prime}}-1}} \lesssim \| v \|_{H^{{s^{_\prime}}-1}}  \lesssim \frac{1}{\ell}.
\end{equation*}
Therefore, using \eqref{brute schatten norm}, we have\footnote{The basis constructed in Proposition \ref{the basis proposition} is not $H^{s^{_\prime}}$-orthonormal, but we may perform a Gram-Schmidt process in each finite dimensional mutual eigenspace of $\curl$ and the Lie bracket by $e_1$. We continue to denote this orthonormal basis as 
before by an abuse of notation.}
\begin{align*}
\norm{K_0(u_0)}_{\mathcal{S}^p}^p = \sum_{\ell=1}^\infty \sum_{v \in \mathcal{E}_{2\ell}^\ell \cup \mathcal{F}_{2\ell}^\ell}\norm{K_0(u_0)v}_{H^{s^{_\prime}}}^p \simeq \sum_{\ell=1}^\infty \ell \frac{1}{\ell^p}
\end{align*}
which implies that $K_0(u_0)$ is Schatten $p$-class on $T_e\Amu^{s^{_\prime}}$ for $p > 2$, just as in 2D ideal hydrodynamics.

\section{The Non-Compact \texorpdfstring{$K_r(u_0)$}{K0}}\label{measures of non-compactness}
As shown in the previous section, in the three (resp. two)-dimensional setting, the operator $K_r$ for the Euler (resp. standard surface quasi-geostrophic) equations, may not belong to any Schatten-von Neumann class. Indeed, it can even fail to be compact, cf. \cite{ebin2006singularities, washabaugh2016the}. The question remains if this phenomenon is in some sense borderline and, if so, can this be leveraged to argue that the failure of Fredholmness for the exponential maps in these settings is of a similar nature.

Here we show that the failure of compactness is essential. We demonstrate this from several perspectives. To this end, we revisit the examples from the proofs of Theorems \ref{2D Example} and \ref{3D Example} which, for the purpose of this section, are fundamentally similar.

Letting $M$ be the round two-sphere, $u_0 = \partial_\theta$ and $r=-\frac{1}{2}$ (the standard SQG equation) we can extract from \eqref{eq_eigen2} and \eqref{rotation eigenvalues} the corresponding operator
\begin{equation*}
    w \mapsto K_{{\scriptscriptstyle -\frac{1}{2}}}(\partial_\theta)w = \sqrt{2} \ \symp \Delta^{-\frac{1}{2}} \partial_\theta \psi_w
\end{equation*}
along with its point spectrum
\begin{equation}\label{SQG eigenvalues}
    \lambda_{\ell}^m = m\left(\sqrt{\frac{2}{\ell(\ell+1)}} \ \right) i, \quad \text{for } \ \ell \geq1 \ , \ -\ell \leq m \leq \ell. 
\end{equation}

Similarly, taking $M$ to be the round three-sphere and $u_0 = e_1$ the Hopf field as our initial data with $r=0$ (the setting of the 3D Euler equations) from \eqref{3D K_0}, \eqref{Hopf field eigenvalues on E_k^m} and \eqref{Hopf field eigenvalues on F_k^m} we obtain the operator
\begin{equation*}
    w \mapsto K_0(e_1)w = -2 \curl^{-1}\big[w, e_1\big]
\end{equation*}
and its point spectrum
\begin{equation}\label{Hopf field eigenvalues}
\begin{aligned}
    \lambda_k^m &= -i (2m-k)\left(\frac{2}{k+2}\right), 
        & \quad \text{for }&\ k \geq 1 ,\ 0 \leq m \leq k \\ 
    \tilde{\lambda}_k^m &= i(2m-k)\left(\frac{2}{k}\right), 
        & \quad \text{for }&\ k \geq 2 ,\ 1 \leq m \leq k-1.
\end{aligned}
\end{equation}
along with the eigenvalues $\pm2i$ associated to eigenvectors from $\mathcal{E}_0^0$ occuring with multiplicity one.

Notice the eigenvalues in \eqref{SQG eigenvalues} and \eqref{Hopf field eigenvalues} densely fill the imaginary intervals $\big[\!-\!\sqrt{2}i, \sqrt{2}i\big]$ and $\big[\!-\!2i, 2i\big]$ respectively, so it is clear that in both examples, $K_r(u_0)$ is not compact. Moreover, in both cases, there exist infinite dimensional subspaces where $K_r(u_0)$ is an isomorphism, i.e., they are not strictly singular\footnote{Recall that an operator is strictly singular if it has no bounded inverse on any infinite dimensional subspace.}.
\subsection{Measures of non-compactness} 
In this subsection we consider the Hausdorff measure of non-compactness\footnote{Equivalent to Kuratowski's measure of non-compactness, as well as other notions.} which, for a bounded linear operator $T$ on a Banach space $X$, is given by
\begin{equation}
    \mu_H(T) = \inf \big\{\rho>0 : T(B_X(0,1)) \ \text{admits a finite cover by balls of radius} \ \rho \big\},
\end{equation}
along with the essential spectral radius
\begin{equation}\label{essential spectral radius}
    r_e(T) = \sup \big\{ \abs{\lambda} : T - \lambda \hspace{0.5mm} \text{Id} \hspace{0.2cm} \text{is not Fredholm} \big\},
\end{equation}
which, heuristically, measures the largest part of the spectrum that persists under compact perturbations. We determine both quantities precisely for the examples above.



\begin{theorem}
    If $M$ is the round two-sphere, $u_0 = \partial_\theta$ and $r=-\frac{1}{2}$ (the standard SQG equation), then
    \[\mu_H\big(K_{{\scriptscriptstyle -\frac{1}{2}}}(\partial_\theta)\big) = r_e\big(K_{{\scriptscriptstyle -\frac{1}{2}}}(\partial_\theta)\big) = \sqrt{2}.\]
    
    If $M$ is the round three-sphere, $u_0 = e_1$ is the Hopf field, and $r=0$ (the 3D Euler equations), then
    \[\mu_H\big(K_{0}(e_1)\big) = r_e\big(K_{0}(e_1)\big) = 2.\]
\end{theorem}

\begin{proof}
    Immediate from the spectra in \eqref{SQG eigenvalues} and \eqref{Hopf field eigenvalues}.
\end{proof}
In fact, the essential spectral radius can be related to any measure of non-compactness $\kappa$ by the formula
$$r_e(T) = \displaystyle \lim_{n\to\infty} \kappa(T^n)^{\frac{1}{n}},$$
indicating that the failure is robust, cf. \cite{nussbaum1970the}.

\subsection{Density of non-decaying eigenvalues} For fixed $\ell \in \Z_{\geq1}$, there are $2\ell+1$ distinct eigenvalues in the spectrum of $K_{{\scriptscriptstyle -\frac{1}{2}}}(\partial_\theta)$. We examine the portion of the spectrum which does not decay as $\ell \to \infty$. In particular, for $\varepsilon>0$ define the density
\begin{equation*}
    \varrho(\ell, \varepsilon)= \dfrac{\# \big\{\lambda_{\ell}^m : \abs{\lambda_\ell^m}>\varepsilon\big\}}{2\ell},
\end{equation*}
where $\#$ denotes the cardinality of the set. From \eqref{SQG eigenvalues}, a simple calculation shows that for any $0<\varepsilon<\sqrt{2}$ we have
\begin{equation*}
    \varrho(\varepsilon) = \lim_{\ell \to \infty} \varrho(\ell, \varepsilon) = 1 - \frac{\varepsilon}{\sqrt{2}}.
\end{equation*}
It follows that $\varrho(\varepsilon) \rightarrow 1$ as $\varepsilon \rightarrow 0$ which means that the set of eigenvalues that do not decay has full density. A similar analysis can be carried out for the spectrum \eqref{Hopf field eigenvalues}, with the same conclusion.

\section{Gaussian Nullity of Singular Sets}\label{Gaussian nullity of singular sets}

We are now ready to present the main results of this paper.  We shall denote by
\begin{equation}
    \mathcal{C}_e^s(M) = \Big\{ u_0 \in T_e\Dex^s(M) \, : \, d\exp_e(u_0) : T_e\Dex^s(M) \rightarrow T_{\exp_e(u_0)}\Dex^s(M) \text{ is not invertible} \Big\}
\end{equation}
the singular sets (sets of conjugate vectors) of the Riemannian exponential maps $\exp_e$ of the right-invariant metrics induced by \eqref{homogeneous sobolev metric} on the group of exact volume-preserving diffeomorphisms.

Our main theorem in the case of ideal hydrodynamics is the following.

\begin{theorem}\label{conjugate theorem}
    For $s>3$ the singular set of the $L^2$ Riemannian exponential map on the group of Sobolev $H^s$ exact volume-preserving diffeomorphisms of the flat two-torus is Gaussian null. 
\end{theorem}
The corresponding result for the family of generalized surface quasi-geostrophic equations is

\begin{theorem}\label{conjugate theorem SQG}
    For $s>3$ and $-\frac{1}{2} < r < 0$ the singular set of the $H^r$ Riemannian exponential map on the group of Sobolev $H^s$ exact volume-preserving diffeomorphisms of the flat two-torus is Gaussian null.
\end{theorem}

The restriction to the flat torus is made for simplicity and is purely an artifact of the tools we employ. We avoid invoking it until necessary. However, the results almost certainly hold for an arbitrary closed surface $M$. We present the proof for the setting of two-dimensional Euler equations ($r=0$). The other settings can be treated in a similar fashion.

\begin{proof}[Proof of Theorem \ref{conjugate theorem}]
Fix $s>3$. By the inverse function theorem, the singular set $\mathcal{C}_e^s$ is closed in the $H^s$ topology and hence Borel. We will show that it is a subset of a Gaussian null set which, as explained in subsection \ref{gaussian null sets}, will imply the result. To this end, consider $u_0 \in T_e\Dex^s$ and let $\gamma(t) = \exp_e(tu_0)$ and $\eta = \gamma(1)$. Note that the smoothness of the exponential map guarantees that $d\exp_e(u_0) : T_e\Dex^s \rightarrow T_{\eta}\Dex^s$ is a bounded linear operator. Furthermore, it is a Fredholm operator of index 0. Hence, $u_0$ is an element of $\mathcal{C}_e^s$ if and only if $d\exp_e(u_0) : T_e\Dex^s \rightarrow T_{\eta}\Dex^s$ is not injective.

Recall the decomposition \eqref{dexp decomposition} from the preliminaries
\begin{equation*}
    d\exp_e(u_0)w = d_eL_{\eta} \big(\Omega(u_0) w + \Gamma(u_0) w\big)
\end{equation*}
and, for $1 < s^{_\prime} \leq s-2$, consider the set
\begin{equation}
    \mathcal{C}_e^{s^{_\prime}}(M) = \Big\{ u_0 \in T_e\Dex^s(M) \, : \, d\exp_e(u_0) : T_e\Dex^{s^{_\prime}}(M) \rightarrow T_{\exp_e(u_0)}\Dex^{s^{_\prime}}(M) \text{ is not invertible} \Big\}. 
\end{equation}

As both $d_eL_\eta : T_e\Dex^{s^{_\prime}} \rightarrow T_\eta\Dex^{s^{_\prime}}$ and $\Omega(u_0) : T_e\Dex^{s^{_\prime}} \rightarrow T_e\Dex^{s^{_\prime}}$ are linear isomorphisms, the derivative
\[d \exp_e (u_0): T_e\Dex^{s^{_\prime}} \rightarrow T_\eta\Dex^{s^{_\prime}}\]
is not injective if and only if the operator
\[\mathrm{Id} + \Omega(u_0)^{-1}\Gamma(u_0) : T_e\Dex^{s^{_\prime}} \rightarrow T_\eta\Dex^{s^{_\prime}}\]
is not injective.

Recall from \eqref{eq_det_p} that for a separable Hilbert space $\mathcal{H}$ and a Schatten-von Neumann $p$-class operator $T$, the map $\mathrm{Id} + T$ is invertible if and only if $\det_p(T) \neq 0$ for any $1 \leq p < \infty$. Furthermore, the determinant is analytic as a map from the Schatten $p$-class equipped with the norm $\norm{\cdot}_{\mathcal{S}_p}$ into $\R$.

\begin{lemma}
    For $u_0 \in T_e\Dex^s(M)$ and $p>2$ the map $\Omega(u_0)^{-1}\Gamma(u_0) \in \mathcal{S}_p\big(T_e\Dex^{s^{_\prime}}(M)\big)$.
\end{lemma}
\begin{proof}
    The fact that $K(u_0) \in \mathcal{S}_p\big(T_e\Dex^{s^{_\prime}}\big)$ for $p>2$ was established in Theorem \ref{main_thm_2d}. The result follows by Lemma \ref{schatten_lemma} applied to 
    \begin{equation*}
        \Gamma(u_0) w = \int_0^1 \Lambda^{-1}(u_0, t) \ K(u_0) \ \Phi(u_0, t) w \ dt
    \end{equation*}
    and the two-sided ideal property \eqref{eq_two_sided_ideal}. 
\end{proof}
An immediate consequence is the following description of the set $\mathcal{C}_e^{s^{_\prime}}$.
\begin{lemma}
    For $u_0 \in T_e\Dex^s(M)$ and $p>2$ the operator $\mathrm{Id} + \Omega(u_0)^{-1}\Gamma(u_0)$ is not injective on $T_e\Dex^{s^{_\prime}}(M)$ if and only if $\mathcal{G}(u_0) = \det_p(\Omega(u_0)^{-1}\Gamma(u_0))=0$.
\end{lemma}
Hence $\mathcal{C}_e^{s^{_\prime}} = \mathcal{G}^{-1}(0)$ and, as outlined in the preliminaries, the proof of Theorem \ref{conjugate theorem} will be complete once we show that the map $\mathcal{G}: T_e\Dex^s \rightarrow \R$ is analytic. To this end, consider the decomposition
\begin{equation}\label{analytic decomp}
    \mathcal{G} = \mathrm{det}_p \circ \Delta^{-\frac{1}{2}} \circ \Delta^{\frac{1}{2}}\Omega(u_0)^{-1}\Gamma(u_0).
\end{equation}
We will treat each factor in \eqref{analytic decomp} separately. We first introduce an analytic chart for $\Dex^s$.

\begin{lemma}\label{analytic lemma}
    There exist open neighbourhoods of zero $\mathcal{V}_e \subset T_e\Dex^s(\T^2)$ and the identity $\mathcal{U}_e \subset \Dex^s(\T^2)$ such that the map
    \begin{equation}\label{analytic chart}
        \Upsilon :  \mathcal{V}_e \rightarrow \mathcal{U}_e, \quad v \mapsto \gamma_v = e + v + \nabla \phi_v
    \end{equation}
    defines a chart around $e \in \Dex^s(\T^2)$. Furthermore, the map $v \mapsto \nabla \phi_v$ is analytic.
\end{lemma}
\begin{proof}
We recall the construction given in   \cite{shnirelman2012on} and \cite{truong2025thesis}. 
As the full diffeomorphism group is an open subset of the model space, the map $w \mapsto e + w$ defines a chart around $e \in \D^s$. By the Hodge decomposition theorem, the vector field $w$ can be written as $w = v + \nabla\phi$ with $v$ divergence-free and $\phi$ a function. Therefore, it suffices to prove that the volume-preserving constraint $\mathrm{det}(D\gamma_v) = 1$ uniquely determines $\nabla\phi$ near the identity as an analytic function of $v$. From \eqref{analytic chart} this constraint can be written in the form
$$
\phi = F(v, \phi) = \Delta^{-1} p(Dv, D^2 \phi),
$$  
where $p$ is a polynomial of degree $2$ in $Dv$ and $D^2\phi$. Hence, $F$ is analytic in both variables. One can show by direct estimates that for sufficiently small $v$, the map $\phi \mapsto F(v, \phi)$ is a contraction, yielding a unique solution $\phi_v$ via the contraction mapping principle. The analytic implicit function theorem then ensures that the map $v \mapsto \nabla \phi_v$ is analytic.
\end{proof}


We now consider the first factor in \eqref{analytic decomp}.

\begin{lemma}
The map $u_0 \mapsto \Delta^{\frac{1}{2}}\Omega(u_0)^{-1}\Gamma(u_0)$ from $T_e\Dex^s(\T^2)$ into $L\big(T_e\Dex^{s^{_\prime}}(\T^2)\big)$ equipped with the operator norm is analytic.
\end{lemma}

\begin{proof}
Let $u_0 \in T_e\Dex^s$ and recall
\begin{equation*}
    \Omega(u_0) : T_e\Dex^{s^{_\prime}} \rightarrow T_e\Dex^{s^{_\prime}}, \quad w \mapsto \int_0^1 \Lambda^{-1}(u_0, t)w \ dt,
\end{equation*}
where the operator $\Lambda(u_0, t): T_e\Dex^{s^{_\prime}} \rightarrow T_e\Dex^{s^{_\prime}}$ can be expressed as
\begin{equation*}
    \Lambda(u_0, t) : w \mapsto P_e D\gamma(t)^{\dagger} D\gamma(t)w,
\end{equation*}
with $P_e$ denoting the $L^2$ projection onto divergence-free fields and $D\gamma(t)^\dagger$ the transpose of the Jacobi matrix of the flow $\gamma(t) = \exp_e(tu_0)$. Recall that $P_e$ is a bounded linear operator in $H^{s^{_\prime}}$ topology and the $L^2$ exponential map on the group of volume-preserving diffeomorphisms of the flat two-torus depends analytically on the initial data, cf. \cite{shnirelman2012on}. Hence, the analytic dependence of $\Lambda(u_0, t)$ on $u_0$ will follow from the analyticity of the map
\begin{equation*}
\gamma \mapsto D\gamma^\dagger D\gamma
\end{equation*}
from $\Dex^s$ to $L\big(T_e\Dex^{s^{_\prime}}, T_e\D^{s^{_\prime}}\big)$ with the usual operator norm which, when written in the local chart from Lemma \ref{analytic lemma}, becomes 
\begin{align*}
v &\mapsto e + Dv + D\big(\nabla\phi_v\big) + Dv^{\dagger} + Dv^{\dagger}Dv + Dv^{\dagger}D\big(\nabla\phi_v\big) + D\big(\nabla\phi_v\big)^{\dagger} + D\big(\nabla\phi_v\big)^{\dagger}Dv + D\big(\nabla\phi_v\big)^{\dagger}D\big(\nabla\phi_v\big).
\end{align*}
This expression is polynomial in $Dv$ and $D\big(\nabla \phi_v\big)$, so its analyticity follows from Lemma \ref{analytic lemma}. Consequently, $u_0 \mapsto \Lambda^{-1}(u_0, t)$ and therefore $u_0 \mapsto \Omega^{-1}(u_0)$ are analytic into $L\big(T_e\Dex^{s^{_\prime}}\big)$.

Moving on to 
\begin{equation*}
    \Gamma(u_0): T_e\Dex^{s^{_\prime}} \rightarrow T_e\Dex^{s^{_\prime}}, \quad w \mapsto \int_0^1 \Lambda^{-1}(u_0, t) \ K(u_0) \ \Phi(u_0, t) w \ dt,
\end{equation*}
note that the map $u_0 \mapsto K(u_0)$ is a bounded linear operator from $T_e\Dex^s$ to $L\big(T_e\Dex^{s^{_\prime}}\big)$ and hence analytic. As for $u_0 \mapsto \Phi(u_0, t)$, observe that $\Phi$ is a solution of the ODE
\begin{equation*}
\partial_t \Phi(u_0, t) = \Lambda(u_0, t)^{-1} + \Lambda(u_0, t)^{-1}K(u_0)\Phi(u_0, t),
\end{equation*}
whose coefficients depend analytically on $u_0$. By the analytic dependence of solutions of ODEs on parameters, we conclude that $\Phi$ is also analytic and the lemma follows.
\end{proof}

\begin{remark}\label{analyticity_rmk}
Analyticity of Riemannian exponential maps in both time and space variables has been discussed in a number of previous works, cf. \cite{shnirelman2012on, truong2025thesis} for example. Furthermore, the results of \cite{kappeler2007analyticity} enables one to prove an analogue of Theorem \ref{conjugate theorem} for one-dimensional models including the Camassa-Holm equations. For the analytic inverse function theorem and dependence on parameters, see \cite{hille1957functional}.
\end{remark}
Regarding the inverse Laplacian, we have
\begin{lemma}
    For $p>2$, the map $A \mapsto \Delta^{-\frac{1}{2}} \circ A$  from  $L\big(T_e\Dex^{s^{_\prime}}(M)\big)$ equipped with the operator norm into $\mathcal{S}_p\big(T_e\Dex^{s^{_\prime}}(M)\big)$ equipped with the Schatten-$p$ norm is analytic.
\end{lemma}

\begin{proof}
    Recalling that for a closed surface $\Delta^{-\frac{1}{2}}$ is a Schatten-$p$ operator on $T_e\Dex^{s^{_\prime}}$ the result follows immediately from the estimate \eqref{eq_two_sided_ideal}.
\end{proof}
Lastly, we restate the result of \cite{boyd2019on} in our setting.
\begin{lemma}
    For $p>2$ the map $\det_p : \Big( \mathcal{S}_p\big(T_e\Dex^{s^{_\prime}}(M)\big),\, \norm{\cdot}_{\mathcal{S}_p}\!\Big) \rightarrow \R$ is analytic.
\end{lemma}
Hence, combining the above lemmas, we conclude that the map $\mathcal{G} = \det_p \circ \Delta^{-\frac{1}{2}} \circ \Delta^{\frac{1}{2}}\Omega(u_0)^{-1}\Gamma(u_0)$ is analytic. This completes the proof of the theorem.    
\end{proof}

\appendix
\section{A Basis of Curl Eigenfields on \texorpdfstring{$\mathbb{S}^3$}{S3}}\label{a basis of curl eigenfields} We construct the Schauder basis for the complexification of the space of divergence-free vector fields on the round three-sphere which is outlined in the proof of Theorem \ref{3D Example}. We recall the statement here.

\begin{proposition}[\textbf{= Prop. 3.8}] Let $e_1 = -y \partial_x + x \partial_y - w \partial_z + z \partial_w$ be the Hopf field on $\mathbb{S}^3$. There exists a Schauder basis of $\C \otimes T_e\Dmu(\mathbb{S}^3)$ consisting of two families of divergence-free vector fields $\mathcal{E}_k^m$ for $k\geq0$ and $0 \leq m \leq k$ and $\mathcal{F}_{k}^m$ for $k \geq 2$ and $1 \leq m \leq k-1$ such that
    \begin{itemize}
    \item The cardinality $\# \mathcal{E}_0^0 = 3$ and, if $v \in \mathcal{E}_0^0$, then $\curl v = 2v$. \\
        \item For $k\geq 1$ the cardinality $\# \mathcal{E}_k^m =
            \begin{cases}
                2(k+1) \quad \text{if } \ m=0, \\
                k+1 \quad \text{if } \ 1\leq m \leq k-1, \\
                2(k+1) \quad \text{if } \ m=k.
            \end{cases}$ \\
        \item For $k\geq1$ if $v \in \mathcal{E}_k^m$ then $\curl v = (k+2)v$ and $[v,e_1] = -i(2m-k)v$. \\
        \item For $k \geq 2$ the cardinality $\# \mathcal{F}_{k}^m = k+1$. \\
        \item For $k \geq 2$ if $v \in \mathcal{F}_{k}^m$ then $\curl v = -kv$ and $[v,e_1] = -i(2m-k)v$.
    \end{itemize}
\end{proposition}
\begin{proof}

We begin the construction of the basis at the level of functions. It is well-known that the Laplacian on $\mathbb{S}^3$ is symmetric with respect to the $L^2$ metric and that we have the following spectral decomposition
\begin{equation*}
L^2(\mathbb{S}^3) = \mathcal{H}_0 \oplus \mathcal{H}_1 \oplus \mathcal{H}_2 \oplus \cdots 
\end{equation*}
where $\dim \mathcal{H}_k = (k+1)^2$ and $\Delta f = k(k+2)f$ for all $f \in \mathcal{H}_k$. Since the Hopf field $e_1$ is a Killing field, the corresponding directional derivative operator  acting on functions  $\nabla_{e_1}: f \mapsto df(e_1)$  commutes with $\Delta$ and preserves the eigenspaces $\mathcal{H}_k$. Integration by parts shows that $\nabla_{e_1}$ is skew-symmetric, so that for each $k = 0, 1, 2, \ldots$ we have a further decomposition of $\mathcal{H}_k$ into eigenspaces of $\nabla_{e_1}$, necessarily associated with complex eigenvalues.

The group structure on $\mathbb{S}^3$ provides us with a right-invariant moving frame
\begin{equation}\label{real frame}
\begin{aligned}
e_1 &= -y \partial_x + x \partial_y - w \partial_z + z \partial_w, \\
e_2 &= -z\partial_x + w\partial_y + x\partial_z -y \partial_w, \\
e_3 &= -w\partial_x -z\partial_y + y\partial_z + x\partial_w
\end{aligned}
\end{equation}
whose bracket relations are given by
\begin{equation}
    [e_1,e_2] = -2e_3, \quad [e_3,e_1] = -2e_2, \quad [e_2,e_3] = -2e_1.
\end{equation}
Identifying $\mathbb{R}^4 \simeq \mathbb{C}^2$ we introduce complex coordinates
\begin{equation*}
\alpha = x + iy,\qquad \bar{\alpha} = x - iy,\qquad \beta = z + iw,\qquad \bar{\beta} = z - iw
\end{equation*}
in which we express the (positive-definite) Laplacian in $\R^4$
\begin{equation}\label{euclidean laplacian}
    \Delta_{\R^4} = -4 \big(\partial_\alpha \partial_{\bar{\alpha}} + \partial_\beta \partial_{\bar{\beta}}\big)
\end{equation}
and our right-invariant frame
\begin{equation}
\begin{aligned}\label{frame in complex coordinates}
    e_1 &= i \big(\alpha \partial_{\alpha} - \bar{\alpha} \partial_{\bar{\alpha}} + \beta \partial_\beta - \bar{\beta} \partial_{\bar{\beta}}\big) \\
    e_2 &= -\bar{\beta}\partial_\alpha - \beta \partial_{\bar{\alpha}} + \bar{\alpha} \partial_\beta + \alpha \partial_{\bar{\beta}} \\
    e_3 &= i \big(-\bar{\beta}\partial_\alpha + \beta \partial_{\bar{\alpha}} + \bar{\alpha} \partial_\beta - \alpha \partial_{\bar{\beta}}\big).
\end{aligned}
\end{equation}

Consider now the family of homogeneous polynomials in the formal variables $z_1, z_2$
\begin{equation*}
F_{k}^m(\alpha, \bar{\alpha}, \beta, \bar{\beta})(z_1, z_2) = (\alpha z_1 + \beta z_2)^m(- \bar{\beta}z_1 + \bar{\alpha}z_2)^{k-m},
\end{equation*}
for integer $k \geq 0$ and $0 \leq m \leq k$ along with, for $0 \leq j \leq k$, the implicitly defined $Q^m_{kj}(\alpha, \bar{\alpha}, \beta, \bar{\beta})$ given by 
\begin{equation}\label{eq_qkmj_def}
F_k^m(\alpha, \bar{\alpha}, \beta, \bar{\beta})(z_1, z_2) = \sum_{j = 0}^k Q^m_{kj}(\alpha, \bar{\alpha}, \beta, \bar{\beta})z_1^jz_2^{k-j}.
\end{equation}
The following properties are readily verifiable\footnote{Recall that any homogeneous harmonic polynomial on $\R^4$ of order $k$ restricts to an element of $\mathcal{H}_k$.}.
\begin{lemma} For integer $k \geq 0$ and $0 \leq m,j \leq k$ we have
\begin{equation}\label{Q harmonic}
    \Delta_{\R^4}Q_{kj}^m = 0, \quad \Delta_{\mathbb{S}^3}Q_{kj}^m = k(k+2)Q_{kj}^m
\end{equation}
and
\begin{align}\label{Q relations}
    \nabla_{e_1} Q_{kj}^m &= i(2m-k)Q_{kj}^m \nonumber\\
    \nabla_{e_2} Q_{kj}^m &= mQ_{kj}^{m-1} - (k-m)Q_{kj}^{m+1} \\
    \nabla_{e_3} Q_{kj}^m &= imQ_{kj}^{m-1} + i(k-m)Q_{kj}^{m+1} \nonumber
\end{align}
with the convention that $Q_{kj}^{-1} = Q_{kj}^{k+1} = 0$.
\end{lemma}
\begin{proof}
    Using \eqref{euclidean laplacian} and \eqref{frame in complex coordinates}, it is not difficult to show that \eqref{Q harmonic} and \eqref{Q relations} hold when the $Q_{kj}^m$ is replaced with $F_k^m$. From \eqref{eq_qkmj_def} the $Q_{kj}^m$ immediately inherit the same properties.
\end{proof}
Moreover, by induction, one can verify that the $Q^m_{kj}$ are linearly independent. Hence, by a dimension count, we have the following.

\begin{lemma}
Let $E(\lambda)$ denote the eigenspace of the map $f \mapsto \nabla_{e_1}f$ associated with the eigenvalue $\lambda \in \mathbb{C}$. Fix $k \geq 0$ and $0 \leq m \leq k$. Then, the polynomials $Q^m_{kj}$ for $0 \leq j \leq k$ defined by \eqref{eq_qkmj_def} form a basis of
\begin{equation*}
\mathcal{H}_k \cap E\big(i(2m-k)\big).
\end{equation*}
In particular, this space has dimension $k+1$. 
\end{lemma}
\begin{remark}
Note that for fixed $k$ and $m$ the basis $\{ Q^m_{kj} \}_{0 \leq j \leq k}$ need not be  $L^2$-orthogonal. 
\end{remark}

Having completed the necessary groundwork at the level of functions, consider now the right-invariant vector fields
\begin{equation}
    e_1, \quad e_2 - i e_3, \quad e_2 + ie_3
\end{equation}
whose complex linear combinations span the complexification of each tangent space of $\mathbb{S}^3$. These fields satisfy the bracket relations
\begin{equation}\label{complex brackets}
    [e_1, e_2 - ie_3] = -2i (e_2-ie_3), \quad [e_1, e_2 + ie_3] = 2i(e_2+ie_3), \quad [e_2-ie_3,e_2+ie_3]=-4ie_1
\end{equation}
and for $k\geq0$ and $0\leq m,j\leq k$ we have
\begin{equation}\label{more Q relations}
    \nabla_{e_2 - i e_3} Q_{kj}^m = 2mQ_{kj}^{m-1} \quad \text{and} \quad \nabla_{e_2 + i e_3} Q_{kj}^m = -2(k-m)Q_{kj}^{m+1}.
\end{equation}
From these we define
\begin{equation}
    v_{kj1}^m = Q_{kj}^m e_1 \ , \qquad v_{kj2}^m = Q_{kj}^m(e_2-ie_3) \ , \qquad v_{kj3}^m = Q_{kj}^m (e_2 + ie_3).
\end{equation}
for integers $k\geq 0$ and $0 \leq m,j\leq k$. The collection
\begin{equation}
    \bigcup_{k=0}^\infty \ \bigcup_{m,j=0}^k \left\{ v_{kj1}^m , \ v_{kj2}^m, \ v_{kj3}^m   \right\}
\end{equation}
forms a Schauder basis for $\C \otimes T_e\D$ with each grading by $k$ (that is, $\{ v_{kj1}^m , \ v_{kj2}^m, \ v_{kj3}^m \}$) being of (complex) dimension $3(k+1)^2$. Hence, by applying $\curl$, we have that the collection
    \begin{equation}
        \bigcup_{k=0}^\infty \ \bigcup_{m,j=0}^k \big\{ \curl v_{kj1}^m , \ \curl v_{kj2}^m, \ \curl v_{kj3}^m \big\}
    \end{equation}
spans the complexification of the space of divergence-free fields. However, to acquire a basis we must remove the redundancies introduced by gradients. Using the fact that the vectors $e_1, e_2$ and $e_3$ are orthonormal, it is not difficult to show that
\begin{equation}\label{complex basis gradient}
    \nabla f = \big(\nabla_{e_1} f \big) e_1 + \frac{1}{2}\big(\nabla_{e_2 + ie_3} f\big) (e_2 - ie_3) + \frac{1}{2}\big(\nabla_{e_2 - ie_3} f\big) (e_2 +ie_3).
\end{equation}
Using \eqref{Q relations} and \eqref{complex basis gradient} we compute that
\begin{equation}\label{Q gradients}
\begin{aligned}
\nabla Q_{kj}^0 &= -ikv_{kj1}^0 \;-\; kv_{kj2}^1,\\
\nabla Q_{kj}^m &= i(2m-k)v_{kj1}^m - (k-m)v_{kj2}^{m+1} + mv_{kj3}^{m-1}, \quad \text{for } \ 1 \leq m \leq k-1,\\
\nabla Q_{kj}^k &= ikv_{kj1}^k + kv_{kj3}^{k-1}.
\end{aligned}
\end{equation}
As the $\curl$ operator annihilates gradients, we obtain a Schauder basis for $\C \otimes T_e\Dmu$
\begin{equation}\label{initial divergence-free basis}
    \bigcup_{k=0}^\infty \left( \bigcup_{j=0}^k \left( \bigcup\limits_{m=0}^k \left\{\curl v_{kj1}^m\right\} \cup \left\{\curl v_{kj2}^0\right\} \bigcup_{m=2}^k \left\{\curl v_{kj2}^m\right\} \cup \left\{\curl v_{kj3}^k \right\}\right) \right)
\end{equation}
with each grading by $k$ of dimension $2(k+1)^2$.

The next phase of the construction is to adjust this into a basis of eigenfields for the operator $\curl$. For $w = f_1e_1 + f_2e_2 + f_3e_3$ one can compute that
\begin{equation}
    \curl w = \big(\nabla_{e_2}f_3 - \nabla_{e_3}f_2 + 2f_1\big)e_1 + \big(\nabla_{e_3}f_1 - \nabla_{e_1}f_3 + 2f_2\big)e_2 + \big(\nabla_{e_1}f_2 - \nabla_{e_2}f_1 + 2f_3\big)e_3.
\end{equation}
From here it is not difficult to check that for $k \geq 0$ and $0 \leq j \leq k$ we have
\begin{align}\label{V curls}
\curl v_{kj1}^m &= 2v_{kj1}^m + i(k-m)v_{kj2}^{m+1} + i m v_{kj3}^{m-1} \quad &\text{for } \ &1 \leq m \leq k-1, \notag \\
\curl v_{kj2}^{m+1} &= -2i(m+1)v_{kj1}^m + (k-2m)v_{kj2}^{m+1} \quad &\text{for } \ &1 \leq m+1 \leq k, \\
\curl v_{kj3}^{m-1} &= -2i\big(k-m+1)v_{kj1}^m - (k-2m)v_{kj3}^{m-1} \quad &\text{for } \ &0 \leq m-1 \leq k-1 \notag
\end{align}
along with the edge cases
\begin{equation}\label{edge V curls}
\begin{aligned}
\curl v_{kj1}^0 &= 2v_{kj1}^0 + ikv_{kj2}^{1},\\
\curl v_{kj1}^k &= 2v_{kj1}^k + ikv_{kj3}^{k-1},\\
\curl v_{kj2}^{0} &= (k+2)v_{kj2}^{0},\\
\curl v_{kj3}^{k} &= (k+2)v_{kj3}^{k}.
\end{aligned}
\end{equation}
In fact, observe that
\begin{equation*}
    \curl^2 v_{kj1}^0 = (k+2)\curl v_{kj1}^0 \quad \text{and} \quad \curl^2 v_{kj1}^k = (k+2)\curl v_{kj1}^k.
\end{equation*}
Hence, the within the grading by $k$ in \eqref{initial divergence-free basis} the elements
\begin{equation}\label{basis edge sets}
    \bigcup_{j=0}^k\left\{\curl v_{kj1}^0\right\} \cup \left\{\curl v_{kj1}^k\right\} \cup \left\{\curl v_{kj2}^0\right\} \cup \left\{\curl v_{kj3}^k \right\}
\end{equation}
are, for $k\geq1$, already a set of $4(k+1)$ linearly independent $\curl$ eigenfields with eigenvalue $k+2$. To address the other elements within each grading by $k \geq 2$ in \eqref{initial divergence-free basis}
\begin{equation}\label{basis needing decomp}
    \bigcup_{j=0}^k\bigcup_{m=1}^{k-1} \big\{\curl v_{kj1}^m\big\} \cup \big\{\curl v_{kj2}^{m+1}\big\}
\end{equation}
we consider the decomposition\footnote{This decomposition was outlined in an unpublished work of Jason Cantarella (personal correspondence).}
\begin{equation}\label{cantarella decomp}
    (2k+2)w = \big(kw + \curl w\big) + \big((k+2)w - \curl w\big)
\end{equation}
which has the property that, for $k\geq0$, if $w= f_1e_1 + f_2e_2 + f_3e_3$ is divergence-free with $f_1, f_2, f_3 \in \mathcal{H}_k$, then (assuming they are non-zero) $kw + \curl w$ and $(k+2)w - \curl w$ are $\curl$ eigenfields with eigenvalues $k+2$ and $-k$ respectively. 

Applying the decomposition \eqref{cantarella decomp} to the basis elements in \eqref{basis needing decomp} and using \eqref{V curls} we have
\begin{equation}\label{V_1 decomp term 1}
    k\curl(v_{kj1}^m) + \curl^2(v_{kj1}^m) = 4(m+1)(k-m+1)v_{kj1}^m + 2i(k-m)(k-m+1)v_{kj2}^{m+1} +2im(m+1)v_{kj3}^{m-1}
\end{equation}
and
\begin{equation}\label{V_1 decomp term 2}
    (k+2)\curl v_{kj1}^m - \curl^2 v_{kj1}^m =  -4m(k-m)v_{kj1}^m + 2im(k-m)v_{kj2}^{m+1} + 2im(k-m)v_{kj3}^{m-1}.
\end{equation}
along with
\begin{equation}\label{V_2 decomp term 1}
    k\curl v_{kj2}^{m+1} + \curl^2 v_{kj2}^{m+1} = -4i(m+1)(k-m+1)v_{kj1}^m + 2(k-m)(k-m+1)v_{kj2}^{m+1} + 2m(m+1)v_{kj3}^{m-1}
\end{equation}
and
\begin{equation}\label{V_2 decomp term 2}
    (k+2)v_{kj2}^{m+1} - \curl^2 v_{kj2}^{m+1} = -4im(m+1)v_{kj1}^m -2m(m+1)v_{kj2}^{m+1} - 2m(m+1)v_{kj3}^{m-1}.
\end{equation}
Notice now that for $1 \leq m \leq k-1$ from \eqref{V_1 decomp term 1} and \eqref{V_2 decomp term 1} we have
\begin{equation}\label{redundancy 1}
    k\curl v_{kj1}^m + \curl^2 v_{kj1}^m = i \Big(k\curl v_{kj2}^{m+1} + \curl^2 v_{kj2}^{m+1}\Big)
\end{equation}
and, similarly, from \eqref{V_1 decomp term 2} and \eqref{V_2 decomp term 2} we have
\begin{equation}\label{redundancy 2}
    (k+2)\curl v_{kj1}^m - \curl^2 v_{kj1}^m = - i\left(\frac{k-m}{m+1}\right)\Big((k+2) \curl v_{kj2}^{m+1} - \curl^2 v_{kj2}^{m+1}\Big).
\end{equation}
Hence, the elements given in \eqref{basis needing decomp} of our basis can be replaced by the $\curl$ eigenfields
\begin{equation}\label{decomposed basis}
    \bigcup_{j=0}^k\bigcup_{m=1}^{k-1} \Big\{k\curl v_{kj1}^m + \curl^2 v_{kj1}^m\Big\} \cup \Big\{(k+2)\curl v_{kj1}^m - \curl^2 v_{kj1}^m\Big\}
\end{equation}
where
\begin{equation*}
    \bigcup_{j=0}^k\bigcup_{m=1}^{k-1} \Big\{k\curl v_{kj1}^m + \curl^2 v_{kj1}^m\Big\}
\end{equation*}
consists of $(k+1)(k-1)$ $\curl$ eigenfields with eigenvalue $k+2$ and
\begin{equation*}
    \bigcup_{j=0}^k\bigcup_{m=1}^{k-1} \Big\{(k+2)\curl v_{kj1}^m - \curl^2 v_{kj1}^m\Big\}
\end{equation*}
consists of $(k+1)(k-1)$ $\curl$ eigenfields with eigenvalue $-k$.

Taking the union over all $k\geq 0$ of \eqref{basis edge sets} and \eqref{decomposed basis} then gives a Schauder basis for $\C \otimes T_e\Dmu$. Notice that there are $(k+3)(k+1)$ $\curl$ eigenfields with eigenvalue $k+2$ for $k\geq0$ and $(k+1)(k-1)$ $\curl$ eigenfields with eigenvalue $-k$ for $k\geq2$.

Lastly, notice that if
\begin{equation}\label{convenient linear combination}
    w = z_1 v_{kj1}^m + z_2 v_{kj2}^{m+1} + z_3 v_{kj3}^{m-1}
\end{equation}
for some $z_1, z_2, z_3 \in \C$, then using \eqref{Q relations}, \eqref{complex brackets} and \eqref{more Q relations} we have that
\begin{equation}
    \big[w,e_1\big] = -i(2m-k)w.
\end{equation}
As every basis element in \eqref{basis edge sets} and \eqref{decomposed basis} has the form \eqref{convenient linear combination} we have a simultaneous eigenbasis for $\curl$ and the Lie bracket $\mathcal{L}_{e_1}: w \mapsto -[w,e_1]$.

Defining now 
\begin{equation}
    \mathcal{E}_0^0 = \big\{e_1, e_2-ie_3, e_2+ie_3\big\}
\end{equation}
and for $k\geq1$
\begin{equation}\label{A basis}
    \mathcal{E}_k^m = \bigcup_{j=0}^k\mathcal{E}_{kj}^m \quad \text{with} \quad \mathcal{E}_{kj}^m =
    \begin{cases}
        \big\{\curl v_{kj1}^0\big\} \cup \big\{\curl v_{kj2}^0\big\} \quad &\text{if } \ m=0 \\
        \big\{k\curl v_{kj1}^m + \curl^2 v_{kj1}^m\big\} \quad &\text{if } \ 1 \leq m \leq k-1 \\
        \big\{\curl v_{kj1}^k\big\} \cup \big\{\curl v_{kj3}^k \big\} \quad &\text{if } \ m=k
    \end{cases}
\end{equation}
along with, for $k\geq2$
\begin{equation}
    \mathcal{F}_{k}^m = \bigcup_{j=0}^k \mathcal{F}_{kj}^m \quad \text{with} \quad \mathcal{F}_{kj}^m =
    \big\{(k+2)\curl v_{kj1}^m - \curl^2 v_{kj1}^m\big\} \quad \text{for } \ 1 \leq m \leq k-1 \\
\end{equation}
yields the desired basis.
\end{proof}

\bibliographystyle{amsplain}
\bibliography{bibliography.bib}

\providecommand{\bysame}{\leavevmode\hbox to3em{\hrulefill}\thinspace}
\providecommand{\MR}{\relax\ifhmode\unskip\space\fi MR }
\providecommand{\MRhref}[2]{%
  \href{http://www.ams.org/mathscinet-getitem?mr=#1}{#2}
}
\providecommand{\href}[2]{#2}
\begin{thebibliography}{10}

\bibitem{arnold1966sur}
V.~Arnold, \emph{Sur la g\'{e}om\'{e}trie diff\'{e}rentielle des groupes de
  {L}ie de dimension infinie et ses applications \`a l'hydrodynamique des
  fluides parfaits}, Ann. Inst. Fourier (Grenoble) \textbf{16} (1966),
  319--361.

\bibitem{arnold2021topological}
V.~Arnold and B.~Khesin, \emph{Topological methods in hydrodynamics}, {S}econd
  ed., Applied Mathematical Sciences, vol. 125, Springer, Cham, 2021.

\bibitem{bauer2024geometric}
M.~Bauer, P.~Heslin, G.~Misio{\l}ek, and S.~Preston, \emph{Geometric analysis
  of the generalized surface quasi-geostrophic equations}, Math. Ann.
  \textbf{390} (2024), no.~3, 4639--4655.

\bibitem{benn2021conjugate}
J.~Benn, \emph{Conjugate points in {$\mathcal{D}^s_\mu (S^2)$}}, J. Geom. Phys.
  \textbf{170} (2021), Paper No. 104369, 14.

\bibitem{bogachev2014onthe}
V.~Bogachev and I.~Malofeev, \emph{On the distribution of smooth functions on
  infinite-dimensional spaces with measures}, Dokl. Akad. Nauk \textbf{454}
  (2014), no.~1, 11--14.

\bibitem{boyd2019on}
C.~Boyd and N.~Snigireva, \emph{On the analyticity of the {F}redholm
  determinant}, Monatsh. Math. \textbf{190} (2019), no.~4, 675--687.

\bibitem{constantin1994formation}
P.~Constantin, A.~Majda, and E.~Tabak, \emph{Formation of strong fronts in the
  {$2$}-{D} quasigeostrophic thermal active scalar}, Nonlinearity \textbf{7}
  (1994), no.~6, 1495--1533.

\bibitem{cordoba2005evidence}
D.~C\'ordoba, M.~Fontelos, A.~Mancho, and J.~Rodrigo, \emph{Evidence of
  singularities for a family of contour dynamicsn equations}, Proc. Natl. Acad.
  Sci. USA \textbf{102} (2005), no.~17, 5949--5952.

\bibitem{cwikel1977weak}
M.~Cwikel, \emph{Weak type estimates for singular values and the number of
  bound states of {S}chr\"odinger operators}, Ann. of Math. (2) \textbf{106}
  (1977), no.~1, 93--100.

\bibitem{drivas2023singularity}
T.~Drivas and T.~Elgindi, \emph{Singularity formation in the incompressible
  {E}uler equation in finite and infinite time}, EMS Surv. Math. Sci.
  \textbf{10} (2023), no.~1, 1--100.

\bibitem{drivas2022conjugate}
T.~Drivas, G.~Misio{\l}ek, B.~Shi, and T.~Yoneda, \emph{Conjugate and cut
  points in ideal fluid motion}, Ann. Math. Qu\'e. \textbf{46} (2022), no.~1,
  207--225.

\bibitem{ebin1970groups}
D.~Ebin and J.~Marsden, \emph{Groups of diffeomorphisms and the motion of an
  incompressible fluid}, Ann. of Math. (2) \textbf{92} (1970), 102--163.

\bibitem{ebin2006singularities}
D.~Ebin, G.~Misio{\l}ek, and S.~Preston, \emph{Singularities of the exponential
  map on the volume-preserving diffeomorphism group}, Geom. Funct. Anal.
  \textbf{16} (2006), no.~4, 850--868.

\bibitem{gohberg1969introduction}
I.~Gohberg and M.~Kreĭn, \emph{Introduction to the theory of linear
  nonselfadjoint operators}, Translations of Mathematical Monographs, vol. Vol.
  18, American Mathematical Society, Providence, RI, 1969, Translated from the
  Russian by A. Feinstein.

\bibitem{hille1957functional}
E.~Hille and R.~Phillips, \emph{Functional analysis and semi-groups}, American
  Mathematical Society Colloquium Publications, vol. Vol. 31, American
  Mathematical Society, Providence, RI, 1957, rev. ed.

\bibitem{holm1998euler}
D.~Holm, , J.~Marsden, and T.~Ratiu, \emph{{E}uler-{P}oincar\'e models of ideal
  fluids with nonlinear dispersion}, Phys. Rev. Lett. \textbf{80} (1998),
  no.~19, 4173--4176.

\bibitem{kappeler2007analyticity}
T.~Kappeler, E.~Loubet, and P.~Topalov, \emph{Analyticity of {R}iemannian
  exponential maps on {D}iff({T})}, J. Lie Theory \textbf{17} (2007), no.~3,
  481--503.

\bibitem{lebrigant2024conjugate}
A.~Le~Brigant and S.~Preston, \emph{Conjugate points along {K}olmogorov flows
  on the torus}, J. Math. Fluid Mech. \textbf{26} (2024), no.~2, Paper No. 24,
  16.

\bibitem{li2022vorticity}
S.~Li, \emph{On analysis of the exponential map of volume-preserving
  diffeomorphism group on closed orientable surfaces through the vorticity},
  arXiv:2204.09497 [math.AP] (2022), 1--20.

\bibitem{lichtenfelz2018normal}
L.~Lichtenfelz, \emph{Normal forms for the {$L^2$} {R}iemannian exponential map
  on diffeomorphism groups}, Int. Math. Res. Not. IMRN (2018), no.~6,
  1730--1753.

\bibitem{lichtenfelz2022axisymmetric}
L.~Lichtenfelz, G.~Misio{\l}ek, and S.~Preston, \emph{Axisymmetric
  diffeomorphisms and ideal fluids on {R}iemannian 3-manifolds}, Int. Math.
  Res. Not. IMRN (2022), no.~1, 446--485.

\bibitem{majda2002vorticity}
A.~Majda and A.~Bertozzi, \emph{Vorticity and incompressible flow}, Cambridge
  Texts in Applied Mathematics, vol.~27, Cambridge University Press, Cambridge,
  2002.

\bibitem{misiolek1996conjugate}
G.~Misio{\l}ek, \emph{Conjugate points in {${\mathscr{D}}_\mu(T^2)$}}, Proc.
  Amer. Math. Soc. \textbf{124} (1996), no.~3, 977--982.

\bibitem{misiolek2015the}
\bysame, \emph{The exponential map near conjugate points in 2{D}
  hydrodynamics}, Arnold Math. J. \textbf{1} (2015), no.~3, 243--251.

\bibitem{misiolek2010fredholm}
G.~Misio{\l}ek and S.~Preston, \emph{Fredholm properties of {R}iemannian
  exponential maps on diffeomorphism groups}, Invent. Math. \textbf{179}
  (2010), no.~1, 191--227.

\bibitem{misiolek2023on}
G.~Misio{\l}ek and T.~Vu, \emph{On continuity properties of solution maps of
  the generalized {SQG} family}, Vietnam J. Math. (2023), 1--9.

\bibitem{nussbaum1970the}
R.~Nussbaum, \emph{The radius of the essential spectrum}, Duke Math. J.
  \textbf{37} (1970), 473--478.

\bibitem{phelps1978gaussian}
R.~Phelps, \emph{Gaussian null sets and differentiability of {L}ipschitz map on
  {B}anach spaces}, Pacific J. Math. \textbf{77} (1978), no.~2, 523--531.

\bibitem{preston2006on}
S.~Preston, \emph{On the volumorphism group, the first conjugate point is
  always the hardest}, Comm. Math. Phys. \textbf{267} (2006), no.~2, 493--513.

\bibitem{rosenberg1997the}
S.~Rosenberg, \emph{The {L}aplacian on a {R}iemannian manifold}, London
  Mathematical Society Student Texts, vol.~31, Cambridge University Press,
  Cambridge, 1997, An introduction to analysis on manifolds.

\bibitem{schatten1960norm}
R.~Schatten, \emph{Norm ideals of completely continuous operators}, Ergebnisse
  der Mathematik und ihrer Grenzgebiete, (N.F.), vol. Heft 27, Springer-Verlag,
  Berlin-G\"ottingen-Heidelberg, 1960.

\bibitem{shnirelman1994generalized}
A.~Shnirelman, \emph{Generalized fluid flows, their approximation and
  applications}, Geom. Funct. Anal. \textbf{4} (1994), no.~5, 586--620.

\bibitem{shnirelman2005microglobal}
\bysame, \emph{Microglobal analysis of the {E}uler equations}, J. Math. Fluid
  Mech. \textbf{7} (2005), S387--S396.

\bibitem{shnirelman2012on}
\bysame, \emph{On the analyticity of particle trajectories in the ideal
  incompressible fluid}, arXiv Preprint (2012), 1--9.

\bibitem{tauchi2022existence}
T.~Tauchi and T.~Yoneda, \emph{Existence of a conjugate point in the
  incompressible {E}uler flow on an ellipsoid}, J. Math. Soc. Japan \textbf{74}
  (2022), no.~2, 629--653.

\bibitem{vizman2008geodesicequations}
C.~Vizman, \emph{Geodesic equations on diffeomorphism groups}, Editura
  Orizonturi Universitare, Timișoara, 2008.

\bibitem{truong2025thesis}
T.~Vu, \emph{Regularity of solution maps of the generalized surface
  quasi-geostrophic equations}, Ph.D. thesis, University of {I}llinois
  {C}hicago, 2025.

\bibitem{washabaugh2016the}
P.~Washabaugh, \emph{The {SQG} equation as a geodesic equation}, Arch. Ration.
  Mech. Anal. \textbf{222} (2016), no.~3, 1269--1284.

\end{thebibliography}

\vfill
\enlargethispage{7\baselineskip}
\end{document}